\documentclass[11pt,a4paper]{amsart}
\usepackage{amsmath}
\usepackage{amssymb}
\usepackage[french,english]{babel}
\usepackage{amsfonts}
\usepackage{mathrsfs}
\usepackage{verbatim}
\usepackage{graphicx}
\usepackage{color}
\newtheorem{theo}{Theorem}

\newtheorem{Lemm}{Lemma}[section]
\newtheorem{Prop}{Proposition}[section]
\newtheorem{Exp}{Example}[section]
\newtheorem{Corol}{Corollary}
\newtheorem{rema}{Remark}
\newenvironment{Rema}{\begin{rema}\normalfont}{\end{rema}}
\newtheorem{defi}{Definition}

\newcommand{\wt}{\widetilde}
\newcommand{\dbtilde}[1]{\overset{\approx}{#1}}

\def\qed{\hfill$\square$\par \vspace{.1cm}}
\newenvironment{Demo}{{\bf Proof.}}{\qed}
\newcommand{\R}{\mathbb{R}}
\newcommand{\C}{\mathbb{C}}
\newcommand{\N}{\mathbb{N}}
\newcommand{\Z}{\mathbb{Z}}

\newcommand{\Bc}{\mathcal{B}}
\newcommand{\Cc}{\mathcal{C}}
\newcommand{\Dc}{\mathcal{D}}
\newcommand{\Ec}{\mathcal{E}}
\newcommand{\Fc}{\mathcal{F}}
\newcommand{\Ic}{\mathcal{I}}

\newcommand{\Kc}{\mathcal{K}}

\newcommand{\Hc}{\mathcal{H}}

\newcommand{\Oc}{\mathcal{O}}

\newcommand{\Sc}{\mathcal{S}}
\newcommand{\Tc}{\mathcal{T}}
\newcommand{\Vc}{\mathcal{V}}

\newcommand{\Ss}{\mathbb{S}}

\newcommand{\Dom}{\mathrm{Dom}}

\newcommand{\seq}[1]{\left<#1\right>}
\newcommand{\eref}[1]{Equation (\ref{#1})}

\usepackage{graphicx}

\newcommand{\e}{{\rm e}}

\usepackage{times}

\begin{document}
\title[Self-adjointness of magnetic Laplacians on triangulations]
{Self-adjointness of magnetic Laplacians on triangulations}

\author[C. Ann\'e]{Colette Ann\'e}
\address{CNRS/Nantes Universit\'e, Laboratoire de Math\'ematique Jean Leray, Facult\'e des Sciences, BP 92208, 44322 Nantes, (France).}

\email{colette.anne@univ-nantes.fr}

\author[H. Ayadi]{Hela Ayadi}

\address{Universit\'e de Monastir, (LR/18ES15)
\& Académie Navale, 7050 Menzel-Bourguiba (Tunisie)}
\email{ halaayadi@yahoo.fr }

\author[Y. Chebbi]{Yassin Chebbi}

\address{Universit\'{e} de Neuch\^{a}tel, Institut de Math\'{e}matiques,
Rue Emile-Argand 11, CH-2000 Neuch\^{a}tel (Switzerland).}
\email{yassin.chebbi@unine.ch}

\author[N. Torki-Hamza]{Nabila Torki-Hamza}
\address{Universit\'e de Monastir, LR/18ES15 \& Institut Sup\'erieur d'Informatique de Mahdia (ISIMa) B.P 05, Campus Universitaire de Mahdia; 5111-Mahdia (Tunisie).}
\email{natorki@gmail.com}


\date{Version of \today}

\begin{abstract}
The notions of magnetic difference operator or magnetic
exterior derivative defined on weighted graphs are discrete analogues of the notion
of covariant derivative on sections of a fibre bundle and its extension on differential
forms. In this paper, we extend these notions to certain 2-simplicial complexes called
triangulations, in a manner compatible with changes of gauge. Then we study the
magnetic Gau\ss-Bonnet operator naturally defined in this context and introduce
the geometric hypothesis of $\chi-$completeness which ensures the essential
self-adjointness of this operator. This gives also the essential self-adjointness
of the magnetic Laplacian on triangulations. Finally we introduce an hypothesis of
bounded curvature for the magnetic potential which permits to caracterize the domain
of the self-adjoint extension.
\end{abstract}


\maketitle



\section{Introduction}
The question of essential self-adjointness of the magnetic Laplace operator was studied in many recent works, such us \cite{Sh}, \cite{CTT}, \cite{GKS}, \cite{M}, \cite{GT} and
\cite{ABDE}. For further references on this topic, see the bibliography of the monograph \cite{KLW}. The first study of the essential self-adjointness of the discrete Laplacian (without magnetic potential) on 1-forms was conducted by the author of \cite{Mas}. Later, in \cite{AT}, the authors gave a new geometric criterion called $\chi-$completeness, which assures the essential self-adjointness of the operator studied in \cite{Mas}. Subsequently, still without magnetic
potential, the author of \cite{Che} generalized the
notion of $\chi-$completeness on weighted {\it triangulations}-- that is,
2-simplicial complexes such as the faces are only triangles. Further aspects of the essential self-adjointness of discrete Laplacian (without magnetic potential) on 1-forms were studied in \cite{BGJ}.
More recently, the authors of \cite{ABDE} introduced the notion
of $\chi_{\alpha}-$completeness related to the magnetic potential $\alpha$,
which is a mix of discrete geometric properties and the behaviour of the
magnetic potential.

The notion of $\chi$-completeness  assures the existence of an exhaustion of the graph
by a family of cut-off functions satisfying certain properties of boundedness
(the definition is given in section \ref{sec:geo}). It
was first introduced in \cite{AT} as a generalization of the existence of an intrinsic  {\it pseudo metric}
with finite balls, a notion introduced in \cite{HKMW}. But it was proved later in \cite[App. A]{LSW}
that these notions are equivalent on locally finite graphs. The notion of intrinsic pseudo metric
was used for instance by \cite{Sch} to obtain self-adjointness in very general settings,
including magnetic operators on graphs, via a Kato type inequality.

In the present work, using the analogy with the smooth case as
presented in \cite{BGV}, we give a generalization of \cite{AT},
\cite{ABDE} and \cite{Che} by introducing a magnetic potential on weighted triangulations
as defined in \cite{Che}. This gives a {\it magnetic triangulation} where we study
self-adjointness of the magnetic Laplace operator in relation with the $\chi$-completeness
and also the geometric meaning on faces of the magnetic field.

To be more precise, on a combinatorial graph, a magnetic potential $\alpha$ is a skewsymmetric
function defined on the edges and with real values.
There are different definitions of the corresponding magnetic diffential $d_\alpha$.
Our guiding principle in this work is to mimic what is done in the smooth case where we can
understand a magnetic potential as a connection on a fiber bundle which defines a covariant
derivative (Definition \ref{def:cov-der}), and the magnetic field corresponds to the curvature
of this connection (Section \ref{magfield}). A good
reference for this point of view is the book \cite{BGV}. In the discrete setting, a fiber
bundle is always trivial, the sections are just functions on the vertices with values in a
given vector space (see  \cite{GT}). In this paper we use complex-valued functions and
the magnetic potential $\alpha$ is understood as defining a parallel transport along
the edges by the quantity $\exp(i\alpha)$.

In the same way we take care of the {\it change of gauge}, so first we define the action
of the gauge group (functions on the vertices with values in $\Ss^1$) and choose a
definition of the magnetic differential $d_\alpha$
with clear equivariant properties with regard to the gauge action, as well as
its formal adjoint $\delta_\alpha$ (see Remark \ref{formal-adjoint} for this notion).
This is done for functions on vertices (0-forms), for skewsymmetric functions on edges
(1-forms) and for skewsymmetric functions on faces (2-forms) (sections \ref{sec:dif-mag},\ref{sec:der-mag}).
Our point of view could be considered as somewhat complicated but we emphasize that it
helps to obtain a coherent framework and good geometric properties, as the notion of
magnetic field.

This magnetic differential defines naturally a Gau\ss-Bonnet operator
$T_\alpha=d_\alpha+\delta_\alpha$ and a magnetic Laplace operator $\Delta_\alpha=T_\alpha^2$.

We show (Theorem \ref{esthm}) that if the triangulation is $\chi-$complete, in the sense of
\cite{Che}, then the operator $T_\alpha$ is essentially self-adjoint, and the magnetic
Laplace operator also (Corollary 7.1). We remark that this result is valid for any
magnetic potential $\alpha$.

But the notion of magnetic potential has a geometric meaning. Using this analogy we
introduce the notion of a magnetic potential with {\it bounded curvature} and apply it to
caracterize the domain of the self-adjoint extension of
the magnetic Gau\ss-Bonnet operator (Theorem \ref{min-bc}).

Finally we give geometric conditions which assure the $\chi$-completeness (Theorem \ref{thh})
and describe examples where it applies, proving that our notion of $\chi$-completeness is
more general than the $\chi_\alpha$-completeness introduced in \cite{ABDE}.

\section{Preliminaries}
\subsection{Basic concepts}
\subsubsection{Graphs}
A graph $\Kc$ is a couple $(\Vc,\Ec)$ where $\Vc$ is a set at most countable whose elements are called \emph{vertices} and $\Ec$, the set of \emph{edges}, is a subset of $\Vc\times\Vc$. It can be considered as a simplicial complex of dimension one.\\ We assume that $\Ec$ is
symmetric and without loops:
\[(x,y)\in\Ec\Longleftrightarrow(y,x)\in\Ec;~~~~x\in\Vc\Longrightarrow(x,x)\notin\Ec.\]
So each edge can be considered with two orientations, and choosing an orientation of the graph
consists of defining a partition of $\Ec$ as follows:
\[\Ec=\Ec^{-}\sqcup\Ec^{+};~~~~(x,y)\in\Ec^{+}\Longleftrightarrow(y,x)\in\Ec^{-}.
\]

Given an edge $e=(x,y),$ we set $e^{-}=x,~e^{+}=y\mbox{ and }-e=(y,x)$
where $e^{-}$ and $e^{+}$ are called the boundary points of $e$. We write $y\sim x$ when
$(x,y)$ or $(y,x)$ is an edge.

A path from $x$ to $y$ is a finite sequence $\gamma=(x_0,....,x_n)$ of vertices in $\Vc$,
such that $x_0=x,~x_n=y$ and $x_{i-1}\sim x_i$ for each $i\in\{1,....,n\}.$
The length of the path $\gamma$ is the number $n.$ If $x_0=x_n,$ we say that
the path is closed or that it is a cycle. If no cycle appears in a path, except maybe the path itself,
the path is called a simple path. A graph is connected if for any two vertices
$x$ and $y,$ there exists a path connecting  $x$ to $y.$
For $x\in\Vc$ we denote $\Vc(x)=\{y\in\Vc;\,(x,y)\in\Ec\}$ the set of its neighbours. A
graph is locally finite if any vertex has a finite number of neighbours.

\begin{defi} A {\it magnetic graph} is a triplet $\Kc_{\alpha}=(\Vc,\Ec,\alpha)$
  where $(\Vc,\Ec)$ defines a graph with set of vertices $\Vc$ and edges $\Ec$ and
  $\alpha$ is a \emph{ magnetic potential} given on the graph: $\alpha$ is a skewsymmetric
  function on $\Ec$ with real values:
  \[\alpha:\Ec\longrightarrow\R\mbox{ satisfies }\forall (x,y)\in\Ec\;\alpha(x,y)=-\alpha(y,x).
  \]
\end{defi}

To simplify the notations, we denote $\alpha(x,y)$ by $\alpha_{xy}.$

{\bf In the sequel, we shall consider all the magnetic graphs as connected, locally finite,
  and without loops.}

\subsubsection{Triangulation}


The notion of a triangulation is defined in \cite{Che}: it is $(\Kc,\Fc)$ where $\Kc=(\Vc,\Ec)$
is a graph and $\Fc\subset\Vc\times\Vc\times\Vc$ a \emph{symmetric set of triangular faces}. By
definition it satisfies
\[ (x,y,z)\in\Fc\Rightarrow (x,y,z) \hbox{ is a cycle of length $3$}, (y,x,z)\in\Fc
\hbox{ and } (y,z,x)\in\Fc.\]
So a cycle of length $3$ is not necessarily a face. A triangulation $\Tc=(\Kc,\Fc)$ can be
considered as a simplicial complex of dimension two. If $f=(x,y,z)\in\Fc$  we say, by an abuse of language,
that  $(y,x,z)$ or  $(y,z,x)$ {\it defines the same face as $f$}, but a face 
has two orientations depending on the signature of the permutation of the vertices.
So $(x,y,z)$ and $(y,z,x)$ define the same orientation, $(x,y,z)$ and $(y,x,z)$ define opposite
orientations. The set of oriented faces  can be seen as a subset of the set of all simple 3-cycles
quotiented by direct permutations. 

\begin{defi}
We say that $\Tc_{\alpha}:=(\Kc_{\alpha},\Fc)$ is a magnetic triangulation,
if $\Kc_{\alpha}$ is a magnetic graph and $\Fc\subset\Vc\times\Vc\times\Vc$
is a symmetric set of triangular faces.
\end{defi}
\begin{rema}We do not need that all the cycles of length 3 define a face. Figure 1 gives an example
  of a triangulation, the white cycles are not faces.
\end{rema}
\begin{figure}[!ht]
\centering
\begin{minipage}[t]{10cm}
\centering
\includegraphics*[width=8cm,height=5cm]{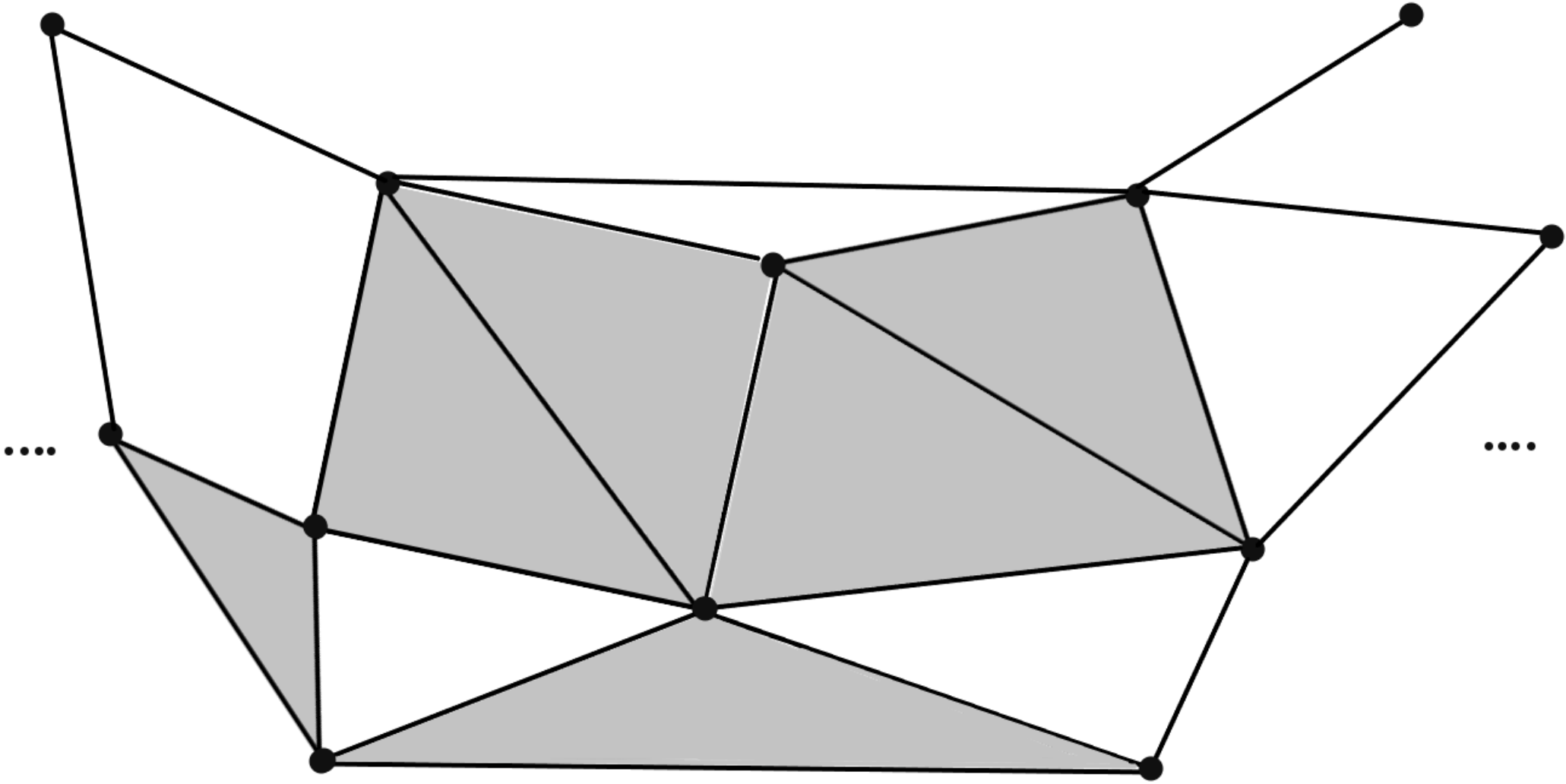}
\caption{A triangulation}
\end{minipage}
\end{figure}
\subsubsection{Weights} To define weighted magnetic triangulations we need weights on a triangulation
$\Tc=(\Kc,\Fc)$, that is
even functions with positive value on vertices, edges and faces:
\begin{itemize}
\item $c: \Vc\rightarrow\R^*_+$ the weight on the vertices.
\item $r:\Ec\rightarrow\R^*_+$ the weight on the edges: $\forall (x,y)\in\Ec,\,r(x,y)=r(y,x) $.
\item $s:\Fc\rightarrow\R^*_+$ the weight on the faces: $\forall (x,y,z)\in\Fc,\,s(x,y,z)=s(x,z,y)=s(y,z,x)$.
\end{itemize}
 
We say that the triangulation $\Tc$ is simple, if, the weights on vertices, edges and faces are all equal to $1.$

We define the set of vertices connected by a face to the edge $(x,y)$ by
$$\Fc_{xy}=\Fc_{yx}:=\{z\in\Vc;~(x,y,z)\in\Fc\}\subseteq\Vc(x)\cap\Vc(y).
$$

The weighted degrees of vertices $deg_{\Vc}$ and of edges $deg_{\Ec}$ are given respectively by:
\[deg_{\Vc}(x):=\dfrac{1}{c(x)}\displaystyle\sum_{y\sim x}r(x,y)\mbox{ and }
deg_{\Ec}(x,y):=\dfrac{1}{r(x,y)}\displaystyle\sum_{z\in\Fc_{xy}}s(x,y,z).
\]

When $\Tc$ is simple, $deg_{\Vc}(x):=deg_{\rm comb}(x)$ the combinatorial
degree, and $deg_{\Ec}(x,y)=|\Fc_{xy}|$, where $|A|=\sharp A$ is the cardinality of the set $A$.

\subsubsection{Holonomy}
In \cite{CTT} and \cite{GT}, the authors introduce the definition of the flux of a magnetic potential:
\begin{defi}Let $\Kc_{\alpha}$ be a magnetic graph, the space of cycles of $\Kc_{\alpha},$ denoted by
  $Z_{1}(\Kc_{\alpha})$ is the $\Z$-module with basis
the geometric cycles $\gamma=(x_0,....,x_n=x_0).$ The holonomy map is the map
$$Hol_{\alpha}:Z_{1}(\Kc_{\alpha})\longrightarrow\R
$$
given by

$$Hol_{\alpha}(\gamma):=\alpha_{x_{0}x_{1}}+\alpha_{x_{1}x_{2}}+........+\alpha_{x_{n-1}x_{n}}.
$$
\end{defi}
If there exists a real function $f$ defined on $\Vc$ such that $\alpha=d_0 f$ where $d_0$
denotes the difference operator: $d_0(f)(x,y)=f(y)-f(x)$, then the holonomy of $\alpha$ is
zero. In the other direction, if the magnetic potential $\alpha$ has no holonomy, then its integration on
paths does not depend on the path joining two points. Because the graph is
connected, this defines a real function whose differential is $\alpha$ (see \cite{GT}).

\begin{defi} We say that the magnetic potential $\alpha$ is {\it trivial} if
  $Hol_{\alpha}=0$.
\end{defi}

In the case of triangulations, any face $(x,y,z)$ defines a cycle $(x,y,z,x)$ so the
holonomy defines a skewsymmetric function $\widehat\alpha$ on $\Fc$:
\begin{equation}
  \label{eq:hol}
\forall (x,y,z)\in\Fc, \;\widehat{\alpha}_{xyz}=Hol_{\alpha}(x,y,z,x)=\alpha_{xy}+\alpha_{yz}+\alpha_{zx}.
\end{equation}

\subsection{Function spaces}
Let $(\Tc_{\alpha},c,r,s)$ be a weighted triangulation.
In this section, we endow Hilbert structures
on the spaces of functions (or cochains) on vertices, edges and faces.

\subsubsection{Hilbert structure on 0-cochains}
The set of 0-cochains is given by:
$$\Cc(\Vc)=\left\{f:\Vc\rightarrow\C\right\},
$$
and its subset of functions of finite support by $\Cc^{c}(\Vc).$

We turn to the Hilbert space:
$$l^2(\Vc):=\left\{f\in\Cc(\Vc);~\displaystyle\sum_{x\in\Vc}c(x)|f(x)|^{2}<\infty\right\},
$$
is endowed with the scalar product given by
$$\langle f_{1},f_{2} \rangle_{l^2(\Vc)}:=\displaystyle\sum_{x\in\Vc}c(x)f_{1}(x)\overline{f_{2}(x)},
$$
for $f_{1},f_{2} \in l^2(\Vc).$
\subsubsection{Hilbert structure on 1-cochains}
The set of complex skewsymmetric 1-cochains is given by:
$$\Cc_{skew}(\Ec)=\left\{\varphi:\Ec\longrightarrow\C;~\varphi(x,y)=-\varphi(y,x)\right\}.
$$

Its subset of functions with finite support is denoted by $\Cc^{c}_{skew}(\Ec)$.
Let us define the Hilbert space
$$l^{2}(\Ec):=\left\{\varphi\in\Cc_{skew}(\Ec);
~\displaystyle\sum_{(x,y)\in\Ec}r(x,y)|\varphi(x,y)|^{2}<\infty\right\}
$$

endowed with the scalar product
$$\langle \varphi_{1},\varphi_{2} \rangle_{l^{2}(\Ec)}
=\dfrac{1}{2}\displaystyle\sum_{(x,y)\in\Ec}r(x,y)\varphi_{1}(x,y)\overline{\varphi_{2}(x,y)}
$$
when $\varphi_{1}$ and $\varphi_{2}$ are in $l^{2}(\Ec).$

\subsubsection{Hilbert structure on 2-cochains}
The set of complex skewsymmetric 2-cochains is given by:
$$\Cc_{skew}(\Fc)=\left\{\psi:\Fc\longrightarrow\C;~\psi(x,y,z)=-\psi(x,z,y)\right\}.
$$

The set of functions with finite support is denoted by $\Cc^{c}_{skew}(\Fc)$.
Let us define the Hilbert space
$$l^{2}(\Fc):=\left\{\psi\in\Cc_{skew}(\Fc);
~\displaystyle\sum_{(x,y,z)\in\Fc}s(x,y,z)~|\psi(x,y,z)|^{2}<\infty\right\}
$$

endowed with the scalar product given by
\[\langle \psi_{1},\psi_{2}\rangle_{l^{2}(\Fc)}
=\dfrac{1}{6}\displaystyle\sum_{(x,y,z)\in\Fc}s(x,y,z)\psi_{1}(x,y,z)\overline{\psi_{2}(x,y,z)}
\]
when $\psi_{1}$ and $\psi_{2}$ are in $l^{2}(\Fc).$

The direct sum of the spaces $l^{2}(\Vc)$, $l^{2}(\Ec)$ and $l^{2}(\Fc)$
can be considered as a Hilbert space denoted by $l^{2}(\Tc_{\alpha})$, indeed
independent on $\alpha$, that is
$$l^{2}(\Tc_{\alpha})=l^{2}(\Vc)\oplus l^{2}(\Ec)\oplus l^{2}(\Fc),
$$

endowed with the scalar product given by
$$\langle (f_{1},\varphi_{1},\psi_{1}),(f_{2},\varphi_{2},\psi_{2}) \rangle_{l^{2}(\Tc_{\alpha})}
=\langle f_{1},f_{2} \rangle_{l^2(\Vc)}
+\langle \varphi_{1},\varphi_{2} \rangle_{l^{2}(\Ec)}
+\langle \psi_{1},\psi_{2} \rangle_{l^{2}(\Fc)}
$$
\begin{Rema}\label{formal-adjoint}The set of functions of finite support $\Cc^c(\Vc)$
  (resp. $\Cc^c(\Ec),\, \Cc^c(\Fc)$)
  are dense in these Hilbert spaces. As a consequence we will take these spaces as the domain
  of the different operators we will introduce later, even if their expression is defined indeed in
  bigger spaces. With this point of view, we can calculate the expression of the adjoint
  of an operator $A$ on these dense subspaces, without worrying about the domain problem.
  The {\bf formal adjoint} of $A$ is the operator $B$ which satisfies, for any functions
  $f,g$ with finite support $\langle A(f),g\rangle=\langle f,B(g)\rangle$, see for instance
  \cite{Bleeker}, or \cite[p.68]{BGV}.
  We recall, see \cite{Reed-Simon}, that a symmetric operator $D$ on a domain $\Dc$ dense in
  the Hilbert space $\Hc$ is essentially self-adjoint if it has one an only one self-adjoint
  extension. If $D$ admits an expression for any element of $\Hc$ then the domain
  of this unique extension is $\{f\in\Hc;\, D(f)\in\Hc\}$.
\end{Rema}

\subsection{A change of gauge} We try here to mimic the approach on manifolds, where one
can consider a magnetic potential as a connection acting on sections of a fiber bundle.
So we consider functions as sections of a $\C$-line bundle on which is defined a connection
$\alpha$. We say that elements of $\Cc(\Vc)$ are {\it sections of the $\alpha$-bundle}.
The {\it gauge group} consists of functions with values in $O(1)=\Ss^1$.
{\it A change of gauge} is then given by any function $f \in \Cc(\mathcal{\Vc},\R)$, which
is considered as acting on sections by multiplication by $\e^{if}$. In the same way elements of
$\Cc_{skew}(\Ec)$ resp. $\Cc_{skew}(\Fc)$ can be considered as 1-forms (resp. 2-forms) with values
in the $\alpha$-bundle (on a manifold $M$ with fiber bundle $E$, $p$-forms of degree $p$ with
value in $E$ are sections of $\Lambda^p(T^\ast M)\otimes E$).
\begin{defi}Let $f \in \Cc(\mathcal{\Vc},\R)$ be a real function on the vertices of the graph.
  It defines the following change of gauge on sections of the $\alpha$-bundle:
  \begin{itemize}
  \item If $g\in \Cc(\mathcal{\Vc})$
\begin{equation}  \label{0form} \e^{if} .g:= \e^{i f}g,\end{equation}
    \item If $\varphi\in \Cc_{skew}(\Ec)$,
\begin{equation}  \label{1form} \e^{if} .\varphi:= \e^{i  \tilde{f}}\varphi,\end{equation}
for the operator of \emph{symmetrization}
\begin{align*}~\widetilde{}: \Cc(\mathcal{V}) &\longrightarrow  \Cc_{sym}(\mathcal{E})\\
  f &\longmapsto\widetilde{f}, \widetilde{f}(e):=\dfrac{f(e^-)+f(e^+)}{2}
\end{align*}
where $\Cc_{sym}(\mathcal{E})=\left\{\varphi:\Ec\longrightarrow\C;~\varphi(x,y)=\varphi(y,x)\right\}.$ 
\item If $\psi\in\Cc_{skew}(\Fc)$,
\begin{equation}  \label{2form} \e^{if}.\psi:= \e^{i   \tilde{\tilde{f}}}\psi,\end{equation}
where the operator
$~\dbtilde{}: \Cc(\mathcal{V}) \longrightarrow  \Cc_{sym}(\mathcal{F})$ is defined
by the formula
$$\dbtilde{f}(x,y,z):=\dfrac{\widetilde{f}(x,y)
   +\widetilde{f}(y,z)+\widetilde{f}(z,x)}{3}=\frac{1}{3}\Big(f(x) +f(y)+f(z)\Big)$$
   and $\Cc_{sym}(\Fc)=\left\{\psi:\Fc\longrightarrow\C;~\psi(x,y,z)=\psi(x,z,y)=\psi(y,z,x)\right\}.$ 
  \end{itemize}
  
\end{defi}
We remark that, on locally finite graphs, these formulae are defined for any function $f\in \Cc(\mathcal{V})$.

 \section{The difference magnetic operator}\label{sec:dif-mag}
  We recall the definition of the difference operator $d^{0}f(e)= f(e^{+})-f(e^{-})$
  which can be defined for any function $f\in \Cc(\mathcal{V})$.

 We look for covariant derivatives, with $\alpha$ as form of connection, which are equivariant
 by change of gauge.  In particular we want that if the magnetic potential is trivial
 ($\alpha=d^{0}f$ for some real function $f\in \Cc(\mathcal{V},\R)$) then a {\it change of gauge}
 permits to come down to the case without magnetic field (see Corollary \ref{trivial}).
 To this way we consider the following definitions, which assure simpler calculus, they are
 different from those of \cite{ABDE} even if they define the same energy form $\|d_{\alpha}^{0}g\|^2$.
 \begin{defi}\label{def:cov-der} On a magnetic triangulation $\Tc_{\alpha}:=(\Kc_{\alpha},\Fc)$
   \textbf{the difference magnetic operator} acting on functions is
   $d_{\alpha}^{0}:\;\Cc(\mathcal{V})\longrightarrow \Cc^{c}_{skew}(\Ec)$ defined by the formula:
   \begin{equation}\label{d0}
     \forall g\in\Cc(\mathcal{V}), (x,y)\in\Ec\quad  d_{\alpha}^{0}g(x,y):=\e^{\frac{i\alpha(y,x)}{2}}g(y)-\e^{\frac{i\alpha(x,y)}{2}}g(x).
      \end{equation}
 \end{defi}
 \begin{Prop}The difference magnetic operator $d_{\alpha}^{0}$ satisfies the
   {\bf gauge invariance}: $\forall f\in \Cc(\Vc,\R)$

   $$\forall g\in \Cc(\mathcal{V})\quad  d_{\alpha+d^0f}^{0}\e^{if}.g =\e^{if}.d_{\alpha}^{0}(g).$$
 \end{Prop}
 \begin{Demo}
   \begin{align*}
     d_{\alpha+d^0f}^{0}(\e^{if}.g)(x,y)
     &=\e^{\frac{i(\alpha(y,x)+f(x)-f(y))}{2}}\e^{if(y)}g(y)-\e^{\frac{i(\alpha(x,y)+f(y)-f(x))}{2}}\e^{if(x)}g(x)\\
     &=\e^{i\wt f(x,y)}\Big(\e^{\frac{i\alpha(y,x)}{2}}g(y)-\e^{\frac{i\alpha(x,y)}{2}} g(x)\Big)\\
     &=\e^{if}. d_{\alpha}^{0}(g)(x,y).
   \end{align*}
 \end{Demo}
 \begin{Corol}\label{trivial}
   Let $\alpha=d^0f$ be a trivial magnetic potential, then the covariant derivative
   $d_{\alpha}^{0}$ is conjugated to the flat one $d^0$:
   \[d_{d^0f}^{0}=\e^{if}. d^0\circ\e^{-if}.
   \]
 \end{Corol}

\begin{defi}\textbf{The co-boundary magnetic operator} is the formal adjoint
    $\delta^0_{\alpha}:\Cc^{c}_{skew}(\Ec)\longrightarrow\Cc^{c}(\Vc)$ of $d^{0}_{\alpha},$
it acts on $\varphi\in\Cc^{c}_{skew}(\Ec)$ by the formula:
$$\forall x\in\Vc,\;(\delta^{0}_{\alpha}\varphi)(x)=\dfrac{1}{c(x)}
\displaystyle\sum_{e,e^{+}=x}r(e)\,{\e}^{i\frac{\alpha_{e}}{2}}\varphi(e).
$$
\end{defi}
Indeed, by definition:
\begin{equation}\label{yc1}\forall f\in \Cc^c(\Vc) ,\,\varphi\in\Cc^{c}_{skew}(\Ec)\quad
\langle d^{0}_{\alpha}f,\varphi\rangle_{l^{2}(\Ec)}
=\langle f,\delta^{0}_{\alpha}\varphi\rangle_{l^{2}(\Vc)}.
\end{equation}
and we calculate
\begin{equation*}
\begin{split}
\langle d^{0}_{\alpha}f,\varphi\rangle_{l^2(\Ec)}
 &=\dfrac{1}{2}\displaystyle\sum_{e\in\Ec}r(e)
\left(\e^{-i\frac{\alpha_{e}}{2}}f(e^{+})-\e^{i\frac{\alpha_{e}}{2}}f(e^{-})\right)\overline{\varphi(e)}\\
& =\sum_{x\in\Vc}f(x)
\left(\sum_{e,e^{+}=x}r(e)\e^{-i\frac{\alpha_{e}}{2}}\overline{\varphi(e)}\right)
=\sum_{x\in\Vc}c(x)f(x)\overline{\left(\dfrac{1}{c(x)}
\sum_{e,e^{+}=x}r(e)\e^{i\frac{\alpha_{e}}{2}}\varphi(e)\right)}.
\end{split}
\end{equation*}

\begin{rema}
  The operators $(d^{0}_{\alpha},\Cc^{c}(\Vc))$ and $(\delta^{0}_{\alpha},\Cc^{c}_{skew}(\Ec))$ are
  closable, as done in \cite{AT}.  It is a simple consequence of the fact that on locally finite
  graphs convergence in norm implies punctual convergence and the operators are \emph{local}.
\end{rema}
Even if it is rather formal in the discrete setting, we consider sections of our $\alpha$-bundle
as a module on scalar functions, and it is important to see what happens for the derivative
under multiplication by a function.
\begin{Prop}Derivation properties\label{devprop0}\\

For $f,g\in\Cc^{c}(\Vc)$ and $\varphi\in\Cc^{c}_{skew}(\Ec)$ it holds
\begin{align}
  \label{Eqd0}
 \forall (x,y)\in\Ec,\quad &d^{0}_{\alpha}(fg)(x,y)=f(y)d^{0}_{\alpha}(g)(x,y)+\e^{i\frac{\alpha(x,y)}{2}}d^{0}(f)(x,y)g(x)
                \Leftrightarrow   \\
 \nonumber\forall e\in\Ec ,\quad \; &d^{0}_{\alpha}(fg)(e)=\wt f(e) d^0_\alpha g(e) +\frac{\e^{i\frac{\alpha(e)}{2}}g(e^-)+\e^{-i\frac{\alpha(e)}{2}}g(e^+)}{2}d^{0}(f)(e)
\\\label{Eqdelta0}\forall x\in\Vc,\quad\;
&\delta^{0}_{\alpha}\left(\wt{f}\varphi\right)(x)
=f(x)\delta^{0}_{\alpha}(\varphi)(x) -\dfrac{1}{2c(x)}\sum_{y\sim x} r(x,y)\e^{i\frac{-\alpha(x,y)}{2}}d^{0}f(x,y)\varphi(x,y).
\end{align}
\end{Prop}
\begin{Demo}

Let $f,g\in\Cc^{c}(\Vc),$ we have that
\begin{equation*}
\begin{split}
d^{0}_{\alpha}(fg)(x,y)& =\e^{i\frac{\alpha(y,x)}{2}}(fg)(y)-\e^{i\frac{\alpha(x,y)}{2}}(fg)(x)\\
& =f(y)\left(\e^{i\frac{\alpha(y,x)}{2}}g(y)-\e^{i\frac{\alpha(x,y)}{2}}g(x)\right)
+\e^{i\frac{\alpha(x,y)}{2}}\left(f(y)g(x)-f(x)g(x)\right)\\
& =f(y)d^{0}_{\alpha}(g)(x,y)+\e^{i\frac{\alpha(x,y)}{2}}d^{0}(f)(x,y)g(x).
\end{split}
\end{equation*}
Furthermore, using the definition $\wt{f}(x,y)=\dfrac{f(x)+f(y)}{2},$ we obtain
\begin{equation*}
 d^{0}_{\alpha}(fg)(x,y)=\wt f(x,y) d^0_\alpha g(x,y) +\frac{\e^{i\frac{\alpha(x,y)}{2}}g(x)+\e^{-i\frac{\alpha(x,y)}{2}}g(y)}{2}d^{0}(f)(x,y).
\end{equation*}
We return to the formula of $\delta^{0}_{\alpha}.$ Given $f\in \Cc^{c}(\Vc)$ and
$\varphi\in\Cc^{c}_{skew}(\Ec),$ as $\wt f(e)=f(e^+)-\frac{1}{2}d^0f(e)$ for any $e\in\Ec$,
we obtain
\begin{equation*}
\begin{split}
\delta^{0}_{\alpha}\left(\wt{f}\varphi\right)(x)
& =\dfrac{1}{c(x)} \sum_{e,e^{+}=x}r(e)\e^{i\frac{\alpha_{e}}{2}}
\left(\widetilde{f}\varphi\right)(e)\\
& =\dfrac{1}{c(x)}f(x)\sum_{y\sim x}r(x,y)\e^{i\frac{\alpha(y,x)}{2}}\varphi(y,x)
+\dfrac{-1}{2c(x)}\sum_{y\sim x}r(x,y)d^0f(y,x)\e^{i\frac{\alpha(y,x)}{2}}\varphi(y,x)\\
& ={f(x)\delta^{0}_{\alpha}(\varphi)(x)}
-\dfrac{1}{2c(x)}\displaystyle\sum_{y\sim x}r(x,y)(d^{0}f)(x,y)\e^{-i\frac{\alpha(x,y)}{2}}\varphi(x,y).
\end{split}
\end{equation*}
\end{Demo}

 \section{The exterior magnetic derivative operator}\label{sec:der-mag}
 The point here is to extend the definition of the covariant derivative with connection form 
 $\alpha$ on 1-forms. We define the \textbf{exterior magnetic derivative operator}
 $d_{\alpha}^{1}:\,\Cc_{skew}(\Ec)\to\Cc_{skew}(\Fc)$ by the formula, if $ \varphi \in \Cc_{skew}(\Ec)$
 and $(x,y,z)\in\Fc$,
 \begin{equation}\label{d-alpha-1}
   d_{\alpha}^{1}(\varphi)(x,y,z):=\e^{\frac{i}{6}(\alpha(x,z)+\alpha(y,z))} \varphi(x,y)+
   \e^{\frac{i}{6}(\alpha(y,x)+\alpha(z,x))} \varphi(y,z)+\e^{\frac{i}{6}(\alpha(z,y)+\alpha(x,y))} \varphi(z,x).
 \end{equation}


It satisfies the following {\bf gauge invariance}. Let $f\in \Cc(\mathcal{V},\R)$
$$d_{\alpha+d^0f}^{1}\e^{if}.\varphi=\e^{if}.d_{\alpha}^{1}(\varphi),\;\;\forall \varphi \in\Cc_{skew}(\Ec).$$
\begin{Demo} Let $(x,y,z)\in\Fc$
  \begin{align*}
    d_{\alpha+d^0f}^{1}(\e^{if}.\varphi)(x,y,z)&= \e^{\frac{i}{6}(\alpha(x,z)+\alpha(y,z)+2f(z)-(f(x)+f(y)))} \e^{\frac{i}{2}(f(x)+f(y))}\varphi(x,y)\\
    &\quad +\e^{\frac{i}{6}(\alpha(y,x)+\alpha(z,x)+2f(x)-(f(z)+f(y)))} \e^{\frac{i}{2}(f(z)+f(y))}\varphi(y,z)\\
    &\quad +\e^{\frac{i}{6}(\alpha(z,y)+\alpha(x,y)+2f(y)-(f(x)+f(z)))} \e^{\frac{i}{2}(f(x)+f(z))}\varphi(z,x)\\
      &=\e^{\frac{i}{3}(f(x)+f(y)+f(z))}d_{\alpha}^{1}(\varphi)(x,y,z).
  \end{align*}
  \end{Demo}
\emph{\textbf{The co-exterior magnetic derivative}} is the formal adjoint of $d^{1}_{\alpha},$
denoted by $\delta^{1}_{\alpha}.$ It satisfies
\begin{align} \label{Eq2.1}
\langle d^{1}_{\alpha}\varphi,\psi   \rangle_{l^{2}(\Fc)}
=\langle \varphi,\delta^{1}_{\alpha}\psi\rangle_{l^{2}(\Ec)},
\end{align}
for all $ (\varphi,\psi)\in\Cc^{c}_{skew}(\Ec)\times\Cc^{c}_{skew}(\Fc).$
\begin{Lemm}\label{lem3.1}
  The formal adjoint $\delta^{1}_{\alpha}:\Cc^{c}_{skew}(\Fc)\longrightarrow\Cc^{c}_{skew}(\Ec)$
  of $d^{1}_{\alpha}$ satisfies, if $\psi\in\Cc^{c}_{skew}(\Fc)$ and $(x,y)\in\Ec$,
\begin{equation*}
  \begin{split}
\delta^{1}_{\alpha}(\psi)(x,y)
& =\dfrac{1}{r(x,y)}\sum_{t\in\Fc_{xy}}s(x,y,t){\rm e}^{\frac{i}{6}(\alpha(t,x)+\alpha(t,y))}\psi(x,y,t).
  \end{split}
\end{equation*}
\end{Lemm}
\begin{Demo}
  Let $\varphi\in\Cc^{c}_{skew}(\Ec)$ and $\psi\in\Cc^{c}_{skew}(\Fc).$ We remark that,
  by \eref{d-alpha-1}, the expression of
  $d^{1}_{\alpha}(\varphi)(x,y,z)$ for $(x,y,z)\in\Fc$ is divided
  into three similar terms. 
  So the equation (\ref{Eq2.1}) gives
  
\begin{equation*}
\begin{split}
\langle d^{1}_{\alpha}\varphi,\psi\rangle_{l^{2}(\Fc)}
& =\dfrac{1}{6}\displaystyle\sum_{(x,y,z)\in\Fc}s(x,y,z)
d^{1}_{\alpha}(\varphi)(x,y,z)\overline{\psi(x,y,z)}\\
& =\dfrac{1}{2}\displaystyle\sum_{(x,y,z)\in\Fc}s(x,y,z)
\e^{\frac{i}{6}(\alpha(x,z)+\alpha(y,z))}\varphi(x,y)\overline{\psi(x,y,z)}\\
& =\dfrac{1}{2}\displaystyle\sum_{(x,y)\in\Ec}r(x,y)\varphi(x,y)
\overline{\left(\dfrac{1}{r(x,y)}\sum_{t\in\Fc_{xy}}s(x,y,t)
\e^{\frac{i}{6}(\alpha(t,x)+\alpha(t,y))}\psi(x,y,t)\right)}.
\end{split}
\end{equation*}
\end{Demo}
\begin{rema}
  The operators $(d^{1}_{\alpha},\Cc^{c}_{skew}(\Ec))$ and $(\delta^{1}_{\alpha},\Cc^{c}_{skew}(\Fc))$ are closable, as done in \cite{Che}.
\end{rema}

To give a derivative property in the same way as Proposition \ref{devprop0}
we need to define the wedge product of a scalar 1-form with a section 1-form.

\begin{defi}
Let $(\xi,\varphi)$ be two 1-forms, we consider $\xi$ as scalar valued and $\varphi$
as a 1-form with values in the $\alpha$-bundle, so we can consider
$\xi\wedge_\alpha \varphi$ as a 2-form with values in the $\alpha$-bundle with the formula
\begin{equation*}
\begin{split}
\left(\xi\wedge_{\alpha}\varphi\right)(x,y,z)
& =\e^{-i(\alpha(z,x)+\alpha(z,y))/6 }\left(\xi(z,x)+\xi(z,y)\right)\varphi(x,y)\\
& +\e^{-i(\alpha(x,y)+\alpha(x,z))/6 }\left(\xi(x,y)+\xi(x,z)\right)\varphi(y,z)\\
& + \e^{-i(\alpha(y,z)+\alpha(y,x))/6 }\left(\xi(y,z)+\xi(y,x)\right)\varphi(z,x).
\end{split}
\end{equation*}
\end{defi}
We remark that for $\alpha=0$, this definition coincides with the wedge product
$\xi\wedge_{disc}\varphi$ given in \cite{Che} for 1-forms with scalar values.

\begin{Prop}\emph{Derivation properties}\label{devprop1}\\

Let $(f,\varphi,\psi)\in\Cc^{c}(\Vc)\times\Cc^{c}_{skew}(\Ec)\times\Cc^{c}_{skew}(\Fc).$ Then we have
\begin{equation}\label{eq:devprop1}
d^1_{\alpha}(\widetilde{f}\varphi)(x,y,z)
=\left({\dbtilde{f}}d^{1}_{\alpha}\varphi\right)(x,y,z)
+\dfrac{1}{6} \left(d^{0}f\wedge_{\alpha}\varphi\right)(x,y,z).
\end{equation}
\begin{equation}\label{eq:devprop2}
\begin{split}
\delta^1_{\alpha}({\dbtilde{f}}\psi)(e)
& =\widetilde{f}(e)(\delta^1_{\alpha}\psi)(e)\\
&\; +\dfrac{1}{6r(e)}\displaystyle\sum_{x\in\Fc_e}s(e,x)\e^{\frac{i}{6}(\alpha(x,e^{-})+\alpha(x,e^{+}))}
  \left(d^0f(e^-,x)+d^0f(e^+,x)\right)\psi(e,x).
\end{split}
\end{equation}
\end{Prop}

\begin{Demo}

\begin{enumerate}
\item [(1)]Let $(f,\varphi)\in\Cc^{c}(\Vc)\times\Cc^{c}_{skew}(\Ec),$ we have
  $\forall (x,y,z)\in\Fc$
\begin{equation*}
\begin{split}
d^1_{\alpha}& (\widetilde{f}\varphi)(x,y,z)-{\dbtilde{f}}(x,y,z)d^1_{\alpha}(\varphi)(x,y,z)=\\
& = \e^{-i(\alpha(z,x)+\alpha(z,y))/6}\left(\frac{f(x)+f(y)}{2}-\frac{f(x)+f(y)+f(z)}{3} \right)\varphi(x,y)\\
  & \; +\e^{-i(\alpha(x,y)+\alpha(x,z))/6 }\left(\frac{f(y)+f(z)}{2}-\frac{f(x)+f(y)+f(z)}{3}
\right)\varphi(y,z)\\
&\; + \e^{-i(\alpha(y,z)+\alpha(y,x))/6 }\left(\frac{f(z)+f(x)}{2}-\frac{f(x)+f(y)+f(z)}{3}
  \right)\varphi(z,x)\\
  & = \e^{-i(\alpha(z,x)+\alpha(z,y))/6}\frac{f(x)+f(y)-2f(z)}{6}
    \varphi(x,y)+
    \e^{-i(\alpha(x,y)+\alpha(x,z))/6 }\frac{f(y)+f(z)-2f(x)}{6}\varphi(y,z)\\
    &\;+ \e^{-i(\alpha(y,z)+\alpha(y,x))/6 }\frac{ f(z)+f(x)-2f(y)}{6}\varphi(z,x)\\
& = \e^{-i(\alpha(z,x)+\alpha(z,y))/6}\dfrac{d^0f(z,x)+d^0f(z,y)}{6}
  \varphi(x,y)+ \e^{-i(\alpha(x,y)+\alpha(x,z))/6 }\dfrac{d^0f(x,y)+d^0f(x,z)}{6}\varphi(y,z)\\
  &\; + \e^{-i(\alpha(y,z)+\alpha(y,x))/6 }\dfrac{ d^0f(y,z)+d^0f(y,x)}{6}\varphi(z,x).
\end{split}
\end{equation*}
  \item [(2)]Let $(f,\psi)\in\Cc^{c}(\Vc)\times\Cc^{c}_{skew}(\Fc).$ In the same way
    as before, using Lemma \ref{lem3.1} one has $\forall e\in\Ec$
  \begin{equation*}
  \begin{split}
  \delta^1_{\alpha}(\dbtilde{f}\psi)(e)-\widetilde{f}(e)\delta^1_{\alpha}(\psi)(e)
  & = \dfrac{1}{r(e)}\displaystyle\sum_{x\in\Fc_e}s(e,x)\e^{\frac{i}{6}(\alpha(x,e^{-})+\alpha(x,e^{+}))}
  \left({\dbtilde{f}}(e,x)-\widetilde{f}(e)\right)\psi(e,x)\\
  & = \dfrac{1}{6r(e)}\displaystyle\sum_{x\in\Fc_e}s(e,x)\e^{\frac{i}{6}(\alpha(x,e^{-})+\alpha(x,e^{+}))}
  \left(d^0f(e^-,x)+d^0f(e^+,x)\right)\psi(e,x).
  \end{split}
  \end{equation*}
\end{enumerate}
\end{Demo}
\section{The magnetic operators}
\subsection{The magnetic Gau\ss -Bonnet operator}
It was originally defined as a square root of the Laplacian. We call magnetic Gau\ss-Bonnet
operator the operator $T_{\alpha}$ defined on
$\Cc^{c}(\Vc)\oplus\Cc^{c}_{skew}(\Ec)\oplus\Cc^{c}_{skew}(\Fc)$ to itself by
\begin{center}
$\begin{pmatrix}
0&\delta^0_{\alpha}&0 \\
d^{0}_{\alpha}&0&\delta^{1}_{\alpha} \\
0&d^{1}_{\alpha}&0
\end{pmatrix}$.
\end{center}
\subsection{The magnetic Laplacian operator}
The magnetic Gau\ss-Bonnet operator $T_{\alpha}$ is of Dirac type and induces the magnetic
Hodge Laplacian on $\Cc^{c}(\Vc)\oplus\Cc^{c}_{skew}(\Ec) \oplus\Cc^{c}_{skew}(\Fc)$ to itself by
$$\Delta_{\alpha}:=T_{\alpha}^{2}.
$$
In general, $\Delta_{\alpha}$ does not preserve the degree of a form, unlike
the usual Hodge Laplacian $(d+\delta)^2$. This default is measured by
the magnetic field.
\subsection{The magnetic field}\label{magfield}
It can be understood as a curvature term which measures how the connection is not flat, it could be
defined as the operator $d^{1}_{\alpha}\circ d^{0}_{\alpha}$ (in the smooth case, a
connection form $A$ defines a covariant derivative written $d_A=d+A$ in local coordinates,
the curvature term is given by $d_A^2=d A+A\wedge A$, see \cite{BGV} Prop. 1.15).
We calculate, for $f\in\Cc^{c}(\Vc)$ and $(x,y,z)\in\Fc$
\begin{multline}
d^{1}_{\alpha}[d^{0}_{\alpha}f](x,y,z)=  \\= \e^{\frac{i}{6}(\alpha(x,z)+\alpha(y,z))}
   d^{0}_{\alpha}f(x,y)+\e^{\frac{i}{6}(\alpha(y,x)+\alpha(z,x))} d^{0}_{\alpha}f(y,z)
   +\e^{\frac{i}{6}(\alpha(z,y)+\alpha(x,y))} d^{0}_{\alpha}f(z,x) \\
     = \e^{\frac{i}{6}(\alpha(x,z)+\alpha(y,z))}
     \left(\e^{\frac{i\alpha(y,x)}{2}}f(y)-\e^{\frac{i\alpha(x,y)}{2}}f(x)\right)\\
     \hspace{2cm}+\e^{\frac{i}{6}(\alpha(y,x)+\alpha(z,x))}
   \left(\e^{\frac{i\alpha(z,y)}{2}}f(z)-\e^{\frac{i\alpha(y,z)}{2}}f(y)\right)\\
     \hspace{4cm}+ \e^{\frac{i}{6}(\alpha(z,y)+\alpha(x,y))}
   \left(\e^{\frac{i\alpha(x,z)}{2}}f(x)-\e^{\frac{i\alpha(z,x)}{2}}f(z)\right) \\
     = -\sin\left(\frac{\widehat{\alpha}(x,y,z)}{6}\right)
     \left(\e^{\frac{i}{3}(\alpha(x,z)+\alpha(x,y))}f(x)
     +\e^{\frac{i}{3}(\alpha(y,z)+\alpha(y,x))}f(y)+
     \e^{\frac{i}{3}(\alpha(z,x)+\alpha(z,y))}f(z)\right)
\end{multline}
where we use to obtain the last line the decomposition $\frac{1}{2}=\frac{1}{6}+\frac{1}{3}$.
So if the holomomy of $\alpha$ is zero, in particular $\widehat{\alpha}=0$ and
$d^{1}_{\alpha}\circ d^{0}_{\alpha}=0$.

Reciprocally if $d^{1}_{\alpha}\circ d^{0}_{\alpha}=0$, by testing the formula on $f$
a Dirac mass at any point $x\in\Vc$ ({\it ie.} $f(x)=1$, and $f(y)=0$ for $y\neq x$) we conclude
that $\sin\left(\dfrac{\widehat\alpha}{6}\right)=0$ on $\Fc$, a kind of
\emph{Bohr-Sommerfeld condition}.

In the same way we can calculate
\begin{equation*}
\begin{split}
   \delta^{0}_{\alpha}[\delta^{1}_{\alpha}\psi](x)
   & = \dfrac{1}{c(x)}\displaystyle\sum_{e,e^{+}=x}r(e)e^{i\frac{\alpha_{e}}{2}}\delta^{1}_{\alpha}\psi(e)\\
     & = \dfrac{1}{c(x)}\displaystyle\sum_{y\sim x}r(x,y)e^{i\frac{\alpha(y,x)}{2}}
     \left(\dfrac{1}{r(x,y)}\sum_{t\in\Fc_{xy}}s(x,y,t)\e^{\frac{i}{6}(\alpha(t,x)+\alpha(t,y))}\psi(y,x,t)\right).
\end{split}
\end{equation*}
We remark that it is less clear (but true) that this term is zero when $\alpha$ has
no holonomy. To see this it is better to remember that
$\delta^{0}_{\alpha}\circ\delta^{1}_{\alpha}$ is the formal adjoint of
$d^{1}_{\alpha}\circ d^{0}_{\alpha}$. It gives another formula, namely
\begin{equation*}
\delta^{0}_{\alpha}\circ\delta^{1}_{\alpha}(\psi)(x)=\dfrac{-3}{c(x)}
\sum_{(x,y,z)\in\Fc} s(x,y,z)\sin\left(\frac{\widehat{\alpha}(x,y,z)}{6}\right)
\e^{\frac{-i}{3}(\alpha(x,z)+\alpha(x,y))}\psi(x,y,z).
\end{equation*}

\begin{rema}\label{gauge} We have proved that, on a magnetic triangulation $\Tc_\alpha$,
  if the magnetic potential has no holonomy, then it is trivial: there exists a
  function $f\in\Cc(\Vc)$ such that $\alpha=d_0f$ and as a consequence by a change of
  gauge the magnetic operators are unitary equivalent to the same operators without
  magnetic potential.

  Moreover, we see now that in this case the magnetic field is nul in the sense that
  $d_\alpha^1\circ d_\alpha^0=0$.
  But the converse is not evident: If the magnetic field is zero, or more strongly if
  $\widehat{\alpha}=0$ we are not sure that the holonomy of $\alpha$ is trivial,
  we need a stronger topological hypothesis as: all the cycles of the graph are
  combinations of boundaries of faces.
\end{rema}


\section{Geometric Hypothesis}\label{sec:geo}
\subsection{Completeness for the magnetic graphs}

We will take the following definition of $\chi-$completeness of a triangulation
\begin{defi}[$\chi-$completeness]\label{de:our}
  Let $\Tc_{\alpha}=(\Kc_{\alpha},\Fc)$ be a weighted magnetic triangulation. $\Tc_{\alpha}$ is
  $\chi$-complete if the underlying triangulation is $\chi$-complete: there exists a positive
  constant $C$ such that the following properties are satisfied:
\begin{itemize}
  \item[($C_1$)] there exists an increasing sequence of finite sets
    $(B_{n})_{n\in\N}$ such that $\Vc=\displaystyle\cup_{n\in\N}B_{n}$
and a sequence of functions $(\chi_{n})_{n\in\N}$ such that
\begin{enumerate}
\item[i)] $\forall n\in\N,\,\chi_{n}\in\Cc^{c}(\Vc),~0\leq\chi_{n}\leq 1$, $\forall  x\in B_n$ $\chi_{n}(x)=1$ ;  
\item[ii)] for all $x\in \Vc$ and $n\in \N$ we have
\[\dfrac{1}{c(x)}\sum_{e\in\Ec,e^{+}=x}r(e)|d^{0}\chi_{n}(e)|^{2}\leq C ;
\]
\end{enumerate}

  \item[($C_2$)] for all $(x,y)\in\Ec$ and $n\in \N$ we have
 \[\dfrac{1}{r(x,y)}\displaystyle\sum_{t\in\Fc_{xy}}s(x,y,t)
  \left|d^{0}\chi_{n}(t,x)+d^{0}\chi_{n}(t,y)\right|^{2}\leq C.
\]
\end{itemize}
\end{defi}

It is indeed the same definition as the one taken in \cite{Che}. We remark that, as a consequence of
$(C_1,\mathbf i)$ one has that for any vertex $x\in\Vc$, $\lim\limits_{n\to\infty}\chi_n(x)=1$.

In \cite{ABDE} is defined the notion of $\chi_{\alpha}-$completeness of a magnetic
weighted graph, a notion mixing geometric properties of the graph with the behavior of
the magnetic potential.

\begin{defi}[$\chi_{\alpha}-$completeness]
The magnetic weighted graph $\Kc_{\alpha}=(\Vc,\Ec,\alpha)$ is $\chi_{\alpha}-$complete
if there exists an increasing sequence of finite sets $(B_{n})_{n\in\N}$ such that $\Vc=\cup_{n\in\N}B_{n}$
and there exist $(\eta_{n})_{n}$ and $(\phi_{n})_n$ with
\begin{enumerate}
\item[i)] $\eta_{n}\in\Cc^{c}(\Vc,\R),~0\leq\eta_{n}\leq 1$ and $\eta_{n}(x)=1$ for all $x\in B_{n}$
\item[ii)] $\phi_{n}\in\Cc(\Vc),$ such that $(\phi_{n})_n$ converges to $0$
\item[iii)]  there exists a non-negative $C$ such that for all $x\in \Vc$ and $n\in \N$, the cut-off
  function $\chi_{n}:=\eta_{n}e^{i\phi_{n}}$ satisfies
\[\dfrac{1}{c(x)}\displaystyle\sum_{e\in\Ec,e^{\pm}=x}r(e)
|d^{0}_{\alpha}\chi_{n}(e)|^{2}\leq C.
\]
\end{enumerate}
\end{defi}
\begin{rema}\label{chialpha}
  We remark that under the hypothesis of $\chi_{\alpha}-$completeness, cut-off functions are allowed to have
complex values. But indeed the graph is $\chi-$complete in the sense of \cite{AT,Che}, using for positive cut-off
functions the absolute values $\eta_n=|\chi_n|$: the triangle inequality gives for any edge $e\in\Ec$
\[|d^0(|\chi_n|)(e)|=||\chi_n(e^+)|-|\chi_n(e^-)||\leq|\chi_n(e^+)-\chi_n(e^-)|=|d^0(\chi_n)(e)|.\]
But the $\chi_{\alpha}-$completeness gives also a constraint on the magnetic potential. Namely, under
the hypothesis of $\chi_{\alpha}-$completeness there
  exists a non-negative constant $C$ such that for all $x\in \Vc$
\begin{equation}
  \label{eq:1}
  \dfrac{1}{c(x)}\displaystyle\sum_{e\in\Ec,e^{\pm}=x}r(e)|\sin(\alpha({e})/2)|^{2}
  \leq C.
\end{equation}
Indeed, let $x \in \Vc$, for $n$ big enough, the vertex $x$ and all its neighbours are in
$B_n$ and one has, because of hypothesis $\mathrm{ii)}$, for any $e\in \Ec, e^\pm=x$:
\begin{equation*}
  \lim_{n\to \infty}d^{0}_{\alpha}\chi_{n}(e)=
  \e^{\frac{-i\alpha(e)}{2}}-\e^{\frac{i\alpha(e)}{2}}=-2i\sin\left(\frac{\alpha(e)}{2}\right).
\end{equation*}
As the graph is locally finite, the conclusion follows. But now, by the general
formula on functions
\begin{equation*}
  d^0_\alpha f(x,y)=\e^{\frac{-i\alpha(x,y)}{2}}f(y)-\e^{\frac{i\alpha(x,y)}{2}}f(x)=
  \e^{\frac{-i\alpha(x,y)}{2}}d^0f(x,y)-2i\sin \left(\frac{\alpha(x,y)}{2}\right)f(x),
\end{equation*}
one has that $\chi_\alpha$-completeness is stronger than $\chi$-completeness in the
sense of Definition \ref{de:our}, indeed it implies $\chi$-completeness and \eref{eq:1}. 
\end{rema}

\subsection{Magnetic graph with bounded curvature}
\begin{defi}
Let $\Tc_{\alpha}=(\Kc_{\alpha},\Fc)$ be a weighted magnetic triangulation. We say that this
magnetic triangulation has a
\emph{bounded curvature} if there exists a
positive constant $C$ such that
\begin{equation}\label{eq:bc}
  \forall x\in\Vc , \;
  \frac{1}{c(x)}\sum_{e\in\Fc_{x}}s(x,e)\sin^2\left(\frac{\widehat{\alpha}(x,e)}{6}\right)\leq C
\end{equation}
\end{defi}
where $\Fc_{x}$ denotes the set of edges connected by a face to the vertex $x.$ For the reason of
this terminology, see Remark \ref{bc} below.


\section{Essential Self-Adjointness}
In \cite{AT} and \cite{Che}, the authors use the $\chi-$completeness hypothesis on a graph
to ensure essential self-adjointness for the Gau\ss -Bonnet operator and the Laplacian.
In this section, with the same idea, we will give a theorem which assures for the
Gau\ss -Bonnet operator $T_\alpha$ to be essentially selfadjoint. We recall the following
definitions for an operator $A$ with domain $\Dc\subset \Hc$ an Hilbert space, the minimal
extension of $A$ is $A_{min}=\overline{A}$ the closure of $A$ and its maximal extension $A_{max}$
has domain $\{f\in\Hc,A(f)\in\Hc\}$ if $A$ has an expression for any function in $\Hc$
(on manifolds it would be the case for differential operators in the sense of distributions,
in our case it is the case as functions on vertices, edges or faces).

Let us begin with this result
\begin{Prop} \label{pro1}
Let $\Tc_\alpha=(\Kc_{\alpha},\Fc)$ be a magnetic $\chi-$complete triangulation, then
the operator $d^0_{\alpha}+\delta^0_{\alpha}$ is essentially self-adjoint on
$\Cc^c(\Vc)\oplus\Cc^c(\Ec).$
\end{Prop}

\begin{Demo}

It suffices to show that $d^0_{\alpha,min}=d^0_{\alpha,max}$ and $\delta^0_{\alpha,min}=\delta^0_{\alpha,max}.$
Indeed, suppose it is proved, $d^0_{\alpha}+\delta^0_{\alpha}$ is a direct sum and if
$F=(f,\varphi)\in \Dom ((d^0_\alpha+\delta^0_\alpha)_{max})$
then $f\in \Dom (d^0_{\alpha,max})$ and $\varphi\in \Dom (\delta^0_{\alpha,max}).$ By
hypothesis, we have then
$f\in \Dom (d^0_{\alpha,min})$ and $\varphi\in \Dom (\delta^0_{\alpha,min}),$ thus $F\in \Dom ((d^0_\alpha+\delta^0_\alpha)_{,min}).$
\begin{enumerate}
\item[1)]Let $f\in \Dom (d^0_{\alpha,max}),$ we will show that
$$\|f-{\chi_{n}}f\|_{l^2(\Vc)}
+\|d^0_{\alpha}\left(f-{\chi_{n}}f\right)\|_{l^2(\Ec)}
\rightarrow 0\mbox{ when }n\rightarrow\infty.
$$
Using the dominated convergence theorem, we have
\begin{align*}\lim_{n\rightarrow \infty}\left( f-{\chi_{n}}f\right)(x) & =\lim_{n\rightarrow \infty}\left( (1-{\chi_{n}})f\right)(x)  =0\\
\hbox{and }\| f-{\chi_{n}}f\|^2_{l^2(\Vc)}
&\leq\displaystyle4\sum_{x\in\Vc}c(x)|f(x)|^2,\;\;f \in l^2(\Vc)
\end{align*}
so $\displaystyle\lim_{n\longrightarrow\infty}
\|f-{\chi_{n}}f\|_{l^2(\Vc)}=0.$
\bigskip

From the derivation formula (\ref{Eqd0}) in Proposition \ref{devprop0}, we have
\begin{equation*}
\begin{split}
d^{0}_{\alpha}((1-{\chi_{n}})f)(x,y)
&=(1-{\chi_{n}})(y)d^{0}_{\alpha}(f)(x,y)+f(x)\e^{\frac{i\alpha(x,y)}{2}}d^{0}((1-{\chi_{n}}))(x,y)\\
&=(1-{\chi_{n}})(y)d^{0}_{\alpha}(f)(x,y)-f(x)\e^{\frac{i\alpha(x,y)}{2}}d^{0}({\chi_{n}})(x,y).
\end{split}
\end{equation*}
We remark that $ d^{0}(1-\chi_{n})(e)=-d^{0}\chi_{n}(e)$ has finite support.
Since $d^0_{\alpha}f\in l^2(\Ec),$ one has by the dominated convergence theorem,
$$\displaystyle\lim_{n\rightarrow\infty}
\|\left(1-{\chi_{n}}\right)d^0_{\alpha}f\|_{l^2(\Ec)}=0.
$$

Now, for the second term, using the triangle inequality, the completeness hypothesis and that the graph is locally finite, we obtain
\begin{eqnarray*}
 \|f(.)d^{0}\chi_{n}\|_{l^2(\Ec)}^2 &= &\frac{1}{2}\sum_{e}r(e) |f(e^+)|^2 | d^{0}\chi_{n}(e)|^2 \\
 &  = & \frac{1}{2}\sum_{x} c(x) |f(x)|^2 \frac{1}{c(x)}\sum_{e,e^{+}=x}r(e) \left( | d^{0}\chi_{n}(e)|^2 \right) \\
  & \leq & C\sum_{x\in\Vc}c(x)|f(x)|^2.
\end{eqnarray*}
But $\forall x\in\Vc$,
\[\lim_{n\to\infty}|f(x)|^2\sum_{e,e^-=x}r(e)|d^{0}\chi_{n}(e)|^2=0.
  \]
We conclude then, by the dominated convergence theorem, that $\displaystyle\lim_{n }\|f(.)d^{0}_{\alpha}(1-\chi_{n})\|_{l^2(\Ec)}=0.$

\item[2)] Let $\varphi\in \Dom (\delta^0_{\alpha,max}),$ we will show that
$$\| {\varphi-\widetilde{\chi_{n}}}\varphi\|_{l^2(\Ec)}
+\|\delta^0_{\alpha}\left(\varphi-\widetilde{\chi_{n}}\varphi\right)\|_{l^2(\Vc)}
\rightarrow 0\mbox{ when }n\rightarrow\infty.
$$
By using dominated convergence theorem, we have
$$ \lim_{n\rightarrow \infty}\left({\varphi-\widetilde{\chi_{n}}}\varphi     \right)(e)  =\lim_{n\rightarrow \infty}\left( (1-\widetilde{\chi_{n}})\varphi\right)(e)  =0$$
$$
\| {\varphi-\widetilde{\chi_{n}}}\varphi\|^2_{l^2(\Ec)}
\leq\displaystyle4\sum_{e\in\Ec}r(e)|\varphi(e)|^2,\;\;\varphi\in l^2(\Ec)
$$
so $\displaystyle\lim_{n\longrightarrow\infty}
\|{\varphi-\widetilde{\chi_{n}}}\varphi\|_{l^2(\Ec)}=0.$

\bigskip

From the derivation formula (\ref{Eqdelta0}) in Proposition \ref{devprop0}, we have
\begin{equation*}
\begin{split}
\delta^{0}_{\alpha}\left((1-\widetilde{\chi_n})\varphi\right)(x)
& =(1-\chi_n)(x)(\delta^{0}_{\alpha}\varphi)(x)
-\dfrac{1}{2c(x)}\sum_{y\sim x}r(x,y)\e^{\frac{-i\alpha(x,y)}{2}}d^{0}(1-\chi_n)(x,y)\varphi(x,y)\\
& =(1-\chi_n)(x)(\delta^{0}_{\alpha}\varphi)(x)
+\dfrac{1}{2c(x)}\sum_{y\sim x}r(x,y)\e^{\frac{-i\alpha(x,y)}{2}}d^{0}(\chi_n)(x,y)\varphi(x,y)
\end{split}
\end{equation*}

On the one hand, we have  $\delta^0_{\alpha}\varphi\in l^2(\Vc),$ then by the dominated
convergence theorem,
$$\displaystyle\lim_{n\rightarrow\infty}
\|\left(1-{\chi_{n}}\right)\delta^0_{\alpha}\varphi\|_{l^2(\Vc)}=0.
$$

On the other hand, fixing $x\in\Vc$, we have
$$A_n(x)=\sum_{e,e^{+}=x}r(e)d^{0}(\chi_n)(e)\varphi(e)\to 0 \hbox{ when } n\to \infty
$$

\begin{eqnarray*}
|A_n(x)|^2 & \leq & \left(\sum_{e,e^{+}=x}r(e)|d^{0}(\chi_n)(e)|^2\right). \left(\sum_{e,e^{+}=x }r(e)|\varphi(e)|^2\right) \\
& \leq & C c(x) \sum_{e,e^{+}=x}r(e)|\varphi(e)|^2.
\end{eqnarray*}

Then, we obtain,
$$c(x)\left( \frac{1}{2 c(x)}\right)^2|A_n(x)|^2
\leq\frac{C}{4}\sum_{e^{+}=x}r(e)|\varphi(e)|^2
$$
which is summable on $\Vc$, as $\varphi\in l^2(\Ec)$.
\end{enumerate}
We conclude then by the dominated convergence theorem that the second term
$\dfrac{1}{2c(x)}A_n(x)$ converges to $0$ in $l^2(\Vc)$.
\end{Demo}
\begin{rema} We see that, because of the good derivative formula given in
  Proposition \ref{devprop0}, the essential selfadjointness of the magnetic
  operator works as if the magnetic potential was zero.
\end{rema}

\begin{Prop} \label{pro2}
Let $\Tc_\alpha=(\Kc_\alpha,\Fc)$ be a $\chi$-complete magnetic triangulation then the operator $d^1_{\alpha}+\delta^1_{\alpha}$
is essentially self-adjoint on $\Cc^c(\Ec)\oplus\Cc^c(\Fc).$
\end{Prop}
\begin{Demo}

Again, as $d^1_{\alpha}+\delta^1_{\alpha}$ is a direct sum, it suffices to show that
$d^1_{\alpha,min}=d^1_{\alpha,max}$ and $\delta^1_{\alpha,min}=\delta^1_{\alpha,max}.$
\begin{enumerate}
\item[1)]Let $\varphi\in \Dom (d^1_{\alpha,max}),$ we will show that
$$\|\varphi-\widetilde{\chi_{n}}\varphi\|_{l^2(\Ec)}
+\|d^1_{\alpha}\left(\varphi-\widetilde{\chi_{n}}\varphi\right)\|_{l^2(\Fc)}
\rightarrow 0\mbox{ when }n\rightarrow\infty.
$$

Using the dominated convergence theorem, we have
\begin{align*}
\forall e\in\Ec,\;&  \lim_{n\rightarrow \infty}\left( \varphi-\widetilde{\chi_{n}}\varphi\right)(e)  =\lim_{n\rightarrow \infty}\left( (1-\widetilde{\chi_{n}})\varphi\right)(e)  =0\\
\| \varphi-&\widetilde{\chi_{n}}\varphi\|^2_{l^2(\Ec)}
\leq\displaystyle4\sum_{e\in\Ec}r(e)|\varphi(e)|^2,\;\;\varphi \in l^2(\Ec)
\end{align*}
so $\displaystyle\lim_{n\longrightarrow\infty}
\|\varphi-\widetilde{\chi_{n}}\varphi\|_{l^2(\Ec)}=0.$
\bigskip

From the derivation formula (\ref{eq:devprop1}) in Proposition \ref{devprop1}, we have
\begin{equation*}
\begin{split}
d^1_{\alpha}\left(\varphi-\widetilde{\chi_{n}}\varphi\right)(x,y,z)
&= d^1_{\alpha}\left(\left(1-\widetilde{\chi_{n}}\right)\varphi\right)(x,y,z)\\
& =\left(\left(1-{\dbtilde{\chi_{n}}}\right)d^{1}_{\alpha}\varphi\right)(x,y,z)
+\dfrac{1}{6}\left(d^{0}(1-\chi_{n})\wedge_{\alpha}\varphi\right)(x,y,z)\\
&=\left(\left(1-{\dbtilde{\chi_{n}}}\right)d^{1}_{\alpha}\varphi\right)(x,y,z)
-\dfrac{1}{6}\left(d^{0}\chi_{n}\wedge_{\alpha}\varphi\right)(x,y,z)
\end{split}
\end{equation*}
using $ d^{0}(1-\chi_{n})(e)=-d^{0}\chi_{n}(e)$.

Since $d^1_{\alpha}\varphi\in l^2(\Fc),$ one has by the dominated convergence theorem,
$$\displaystyle\lim_{n\rightarrow\infty}
\left\|\left(1-{\dbtilde{\chi_{n}}}\right)d^1_{\alpha}\varphi\right\|_{l^2(\Fc)}=0.
$$

For the second term, we recall the definition of $\wedge_\alpha$
\begin{equation*}
\begin{split}
\left(d^{0}\chi_{n}\wedge_{\alpha}\varphi\right)(x,y,z)
& =\e^{-i/6 (\alpha(z,x)+\alpha(z,y))}\left(d^{0}\chi_{n}(z,x)+d^{0}\chi_{n}(z,y)\right)\varphi(x,y)\\
& +\e^{-i/6 (\alpha(x,y)+\alpha(x,z))}\left(d^{0}\chi_{n}(x,y)+d^{0}\chi_{n}(x,z)\right)\varphi(y,z)\\
& +\e^{-i/6 (\alpha(y,z)+\alpha(y,x))}\left(d^{0}\chi_{n}(y,z)+d^{0}\chi_{n}(y,x)\right)\varphi(z,x).
\end{split}
\end{equation*}
We have by property $(C_2)$ of the completeness hypothesis for the triangulation
\begin{equation*}
  \begin{split}
     &\sum_{(e,x)\in\Fc} s(e,x)|\varphi(e)|^2
     \left|d^{0}\chi_{n}(x,e^-)+d^{0}\chi_{n}(x,e^+)\right|^2\\
     & = \sum_{e\in\Ec} r(e)|\varphi(e)|^2 \frac{1}{r(e)}
     \sum_{x\in\Fc_e}s(e,x)\left|d^{0}\chi_{n}(x,e^-)+d^{0}\chi_{n}(x,e^+)\right|^2\\
     & \leq C \sum_{e\in\Ec}r(e)|\varphi(e)|^2.
  \end{split}
\end{equation*}
But for any
$e\in\Ec,\,\lim_{n\to\infty}\sum_{x\in\Fc_e}s(e,x)\left|d^{0}\chi_{n}(x,e^-)+d^{0}\chi_{n}(x,e^+)\right|=0$,
we conclude then again by the dominated convergence theorem that this second term
converges to 0.

\item[2)] Let $\psi\in \Dom (\delta^1_{\alpha, max}),$ we will show that
$$\|\psi-\dbtilde{\chi_{n}}\psi\|_{l^2(\Fc)}
+\|\delta^1_{\alpha}(\psi-\dbtilde{\chi_{n}}\psi)\|_{l^2(\Ec)}
\rightarrow 0\mbox{ when }n\rightarrow\infty.
$$
First, by the dominated convergence theorem, we have
$$ \lim_{n\rightarrow \infty}\left( \psi-\dbtilde{\chi_{n}}\psi\right)
(x,y,z) =
\lim_{n\rightarrow \infty}\left( (1-\dbtilde{\chi_{n}})\psi\right) (x,y,z) =0.$$
\begin{equation*}
\begin{split}
\|\psi-\dbtilde{\chi_{n}}\psi\|^2_{l^2(\Fc)}
& = \frac{1}{6}\displaystyle\sum_{(x,y,z)\in\Fc}s(x,y,z)|1-\dbtilde{\chi_{n}}(x,y,z)|^2
|\psi(x,y,z)|^2\\
& \leq \frac{4}{6} \displaystyle\sum_{(x,y,z)\in\Fc}s(x,y,z)|\psi(x,y,z)|^2=4\|\psi\|^2_{l^2(\Fc)}
.
\end{split}
\end{equation*}
Then,
$$\lim_{n\rightarrow \infty} \|\psi-\dbtilde{\chi_{n}}\psi\|_{l^2(\Fc)}= 0.
$$
Secondly, by the derivation formula (\ref{eq:devprop2}) in Proposition \ref{devprop1}, we have
\begin{equation*}
\begin{split}
M_n(e)& := \delta^1_{\alpha}\left( (1-\dbtilde{\chi_{n}})\psi\right)(e)
-(1-\widetilde{\chi_{n}})(e)(\delta^1_{\alpha}\psi)(e)\\
& = \dfrac{1}{6r(e)}\sum_{x\in\Fc_e}s(e,x)\e^{\frac{i}{6}(\alpha(x,e^{-})+\alpha(x,e^{+}))}
\left(d^0(1-{\chi_{n}})(e^-,x)+d^0(1-{\chi_{n}})(e^+,x)\right)\psi(e,x)\\
& = \dfrac{1}{6r(e)}\sum_{x\in\Fc_e}s(e,x)\e^{\frac{i}{6}(\alpha(x,e^{-})+\alpha(x,e^{+}))}
\left(d^0{\chi_{n}}(x,e^-)+d^0{\chi_{n}}(x,e^+)\right)\psi(e,x)
\end{split}
\end{equation*}
Then, by the hypothesis of completeness,
\begin{equation}
  \label{eq:3}
  \forall e\in \Ec,\;\lim_{n\to\infty}M_n(e)=0.
\end{equation}
On the other hand, using the Cauchy-Schwarz inequality and again the property
$(C_2)$ of $\chi-$completeness, we have for any $e\in\Ec,$
\begin{equation}\label{delta1}
  \begin{split}
    |M_n(e)|^2 &=\left(\frac{1}{6r(e)}\right)^2\left|
     \sum_{x\in \Fc_{e}}s(e,x)\e^{\frac{i}{6}(\alpha(x,e^{-})+\alpha(x,e^{+}))} \left(d^0(\chi_{n})(x,e^-)
     +d^0(\chi_{n})(x,e^+)\right)\psi(e,x)\right|^2 \\
       & \leq \left(\frac{1}{6r(e)}\right)^2\sum_{x\in \Fc_{e}}s(e,x)
       \left|d^0(\chi_{n})(x,e^-)+d^0(\chi_{n})(x,e^+)\right|^2
       \sum_{x\in\Fc_{e}}s(e,x)|\psi(e,x)|^2\\
       & \leq \frac{C}{r(e)} \sum_{x\in\Fc_{e}}s(e,x)|\psi(e,x)|^2.
  \end{split}
\end{equation}
Therefore, by the dominated convergence theorem
\begin{equation*}
\forall n\in \N,\,
 \sum_{e\in \Ec}r(e)|M_n(e)|^2\leq C\|\psi\|^2\Rightarrow
      \lim_{n\to\infty}\sum_{e\in \Ec}r(e)|M_n(e)|^2=0
    \end{equation*}
    Finally, we have by dominated convergence theorem as $\delta^1_{\alpha}\psi\in l^2(\Ec)$
$$\displaystyle\lim_{n\rightarrow\infty}\|\left(1-\widetilde{\chi_{n}}\right)
\delta^1_{\alpha}(\psi)\|=0.
$$
\end{enumerate}
This completes the proof.
\end{Demo}
\begin{theo}\label{esthm}
Let $\Tc_{\alpha}=(\Kc_{\alpha},\Fc)$ be a $\chi-$complete magnetic triangulation then the operator
$T_{\alpha}$ is essentially self-adjoint on $\Cc^c(\Vc)\oplus\Cc^c(\Ec)\oplus\Cc^c(\Fc).$
\end{theo}
\begin{Demo} To show that $T_{\alpha}$ is essentially self-adjoint, we will prove
  that $\Dom (T_{\alpha, max})\subset \Dom (T_{\alpha, min}).$
Let us take the sequence of cut-off functions $\left(\chi_n\right)_n\subseteq
\Cc^c(\Vc)$ assured by the hypothesis of $\chi-$completness of $\Tc_{\alpha}$.
Let $F=(f,\varphi,\psi)\in \Dom (T_{\alpha,max})$, then $F$ and $T_{\alpha}F$ are in $l^2(\Tc_{\alpha}).$ This implies that
$\delta^0_{\alpha}\varphi\in l^2(\Vc),$ $d^0_{\alpha}f+\delta^1_{\alpha}\psi\in l^2(\Ec)$ and $d^1_{\alpha}\varphi\in l^2(\Fc).$
Consequently, by the definition of $\delta^{0}_{\alpha,max}$ and $d^{1}_{\alpha,max}$
we have $\varphi\in \Dom (\delta^{0}_{\alpha,max})\cap \Dom (d^{1}_{\alpha,max}).$ But
the proofs of Proposition \ref{pro1} and Proposition \ref{pro2} show that
$\varphi_n=\wt\chi_n \varphi$ satisfies
\begin{equation*}
  \forall n\in\N\; \varphi_n\in \Cc(\Ec) \hbox{ and
  }\lim_{n\to\infty}\|\varphi-\varphi_n\|^2+\|d^1_\alpha(\varphi-\varphi_n)\|^2
  +\|\delta^0_\alpha(\varphi-\varphi_n)\|^2=0,
\end{equation*}
so $\varphi\in \Dom (\delta^{0}_{\alpha,min})\cap \Dom (d^{1}_{\alpha,min}).$

Now, it remains to prove that $d^0_{\alpha}(\chi_n f)+\delta^1_{\alpha}(\dbtilde{\chi_{n}}\psi)$
converges in $l^2(\Ec)$ to $d^0_\alpha f+\delta^1_\alpha \psi$. Indeed, we need some
derivation formula taken in Proposition
\ref{devprop0} and \eref{delta1} in the proof of Proposition (\ref{pro2}). It gives:
\begin{equation*}
  \begin{split}
     d^0_{\alpha}((1-\chi_n)f)(e)
     & ={(1-\wt{\chi_n})}(e)d^{0}_{\alpha}(f)(e)-\underbrace{\left(\frac{e^{i\frac{\alpha(e)}{2}}f(e^-)+
           e^{-i\frac{\alpha(e)}{2}}f(e^+)}{2}\right)d^0\chi_n(e)}_{\Sc_n(e)}\\
     \delta^1_{\alpha}\left( (1-\dbtilde{\chi_{n}})\psi\right)(e)
     & =(1-\widetilde{\chi_{n}})(e)(\delta^1_{\alpha}\psi)(e)\\
& -\underbrace{\dfrac{1}{6r(e)}\sum_{x\in\Fc_e}s(e,x)\e^{\frac{i}{6}(\alpha(x,e^{-})+\alpha(x,e^{+}))}
\left(d^0({\chi_{n}})(e^-,x)+d^0({\chi_{n}})(e^+,x)\right)\psi(e,x) }_{\Ic_n(e)}.
\end{split}
\end{equation*}

Therefore, we have by the triangle inequality
\begin{equation*}
\begin{split}
&\|d^0_{\alpha}(f-\chi_nf)+\delta^1_{\alpha}(\psi-\dbtilde{\chi_{n}}\psi)\|^2_{l^2(\Ec)}\\
&=\|(1-\widetilde{\chi_n})(d^0_{\alpha}f+\delta^1_{\alpha}\psi)-\Sc_n-\Ic_n\|^2_{l^2(\Ec)}\\
&\leq 3\left(\|(1-\widetilde{\chi_n})(d^0_{\alpha}f+\delta^1_{\alpha}\psi)\|^2_{l^2(\Ec)}
+\|\Sc_n\|^2_{l^2(\Ec)}+\|\Ic_n\|^2_{l^2(\Ec)}\right)
\end{split}
\end{equation*}

Because $d^0_{\alpha}f+\delta^1_{\alpha}\psi\in l^2(\Ec),$ we have
$$\displaystyle\lim_{n\rightarrow\infty}
\|(1-\widetilde{\chi_n})(d^0_{\alpha}f+\delta^1_{\alpha}\psi)\|^2_{l^2(\Ec)}=0.
$$

Because $\psi\in l^2(\Fc)$ we have, as in the proof of  Proposition \ref{pro2}
$$\displaystyle\lim_{n\rightarrow\infty}\|\Ic_n\|^2_{l^2(\Ec)}=0.
$$

Because $f\in l^2(\Vc)$ we have, as in the proof of Proposition \ref{pro1},
$$\displaystyle\lim_{n\rightarrow\infty}\|\Sc_n\|^2_{l^2(\Ec)}=0.
$$
\end{Demo}
\begin{Corol}
Let $\Tc_{\alpha}=(\Kc_{\alpha},\Fc)$ be a $\chi-$complete magnetic triangulation
then the magnetic Laplace operator
$\Delta_{\alpha}$ is essentially self-adjoint on $\Cc^c(\Vc)\oplus\Cc^c(\Ec)\oplus\Cc^c(\Fc).$
\end{Corol}
Knowing that $T_\alpha$ is essentially self-adjoint, we conclude then that
$\Delta_\alpha=T_\alpha^2$ is also essentially self-adjoint in the same way as in
Proposition 13 of \cite{AT}.

\begin{theo}\label{min-bc}
Let $\Tc_{\alpha}=(\Kc_{\alpha},\Fc)$ be a $\chi-$complete triangulation with bounded
curvature
({\it ie.} the property (\ref{eq:bc}) is satisfied) then the operator $T_{\alpha}$ satisfies
$$\Dom (T_{\alpha,min})=\Dom (d^{0}_{\alpha,min})\oplus\left(\Dom (\delta^{0}_{\alpha,min})\cap \Dom (d^{1}_{\alpha,min})\right)
\oplus \Dom (\delta^{1}_{\alpha,min}).
$$\end{theo}
\begin{rema}\label{bc}This theorem is the reason why we call our hypothesis (\ref{eq:bc}) \emph{
    bounded curvature}. Indeed it sounds like, in the smooth case, the theorem which
  says that on a complete manifold with bounded geometry (\emph{ie.} with positive injectivity
  radius and Ricci curvature bounded from below) the Laplace Beltrami
  operator is essentially self adjoint and the domain of its selfadjoint extension is
  the second Sobolev space (the closure of the smooth function with compact support
  for the Sobolev norm of degree 2), see \cite{Hebey} Prop. 2.10.
\end{rema}
\begin{Demo}

\emph{First inclusion.} Let $F=(f,\varphi,\psi)\in \Dom (d^{0}_{\alpha,min})\oplus
\left(\Dom (\delta^{0}_{\alpha,min})\cap \Dom (d^{1}_{\alpha,min})\right)\oplus \Dom (\delta^{1}_{\alpha,min}).$
Then by Proposition \ref{pro1} and Proposition \ref{pro2} we have $(f_n)_n\subseteq\Cc^c(\Vc)$ and
$(\psi_n)_n\subseteq\Cc^c(\Fc)$  such that :
\begin{enumerate}
\item[-]$f_n={\chi_{n}}f\rightarrow f\mbox{ in }l^2(\Vc)$ and
$d^0_{\alpha} f_n \rightarrow d^0_{\alpha,min} f\mbox{ in }l^2(\Ec).$
\item[-]$\psi_n =\dbtilde{\chi_{n}}\psi\rightarrow\psi\mbox{ in }l^2(\Fc)$ and
$\delta^1_{\alpha} \psi_n\rightarrow\delta^1_{\alpha,min}\psi\mbox{ in }l^2(\Ec)$
\end{enumerate}
We have also $\varphi\in \Dom (\delta^{0}_{\alpha,min})\cap \Dom (d^{1}_{\alpha,min})$ so by these two propositions
$(\varphi_{n} =\widetilde{\chi_{n}}\varphi)_n\subseteq\Cc^c(\Ec)$
satisfies
$$\|\varphi-\widetilde{\chi_{n}}\varphi\|_{l^2(\Ec)}
+\|\delta^0_{\alpha}(\varphi-\widetilde{\chi_{n}}\varphi)\|_{l^2(\Vc)}
+\|d^1_{\alpha}(\varphi-\widetilde{\chi_{n}}\varphi)\|_{l^2(\Fc)}\rightarrow 0,
\mbox{ when } n\rightarrow\infty.
$$
Hence, $F_n=(f_n,\varphi_{n},\psi_n)$ satisfies
$$F_n\rightarrow F\mbox{ in }l^2(\Tc_{\alpha}),~T_{\alpha}F_n\rightarrow T_{\alpha, min}F
\mbox{ in }l^2(\Tc_{\alpha}),
$$
where $T_{\alpha,min}F=(\delta^0_{\alpha,min} \varphi,d^0_{\alpha,min}f
+\delta^1_{\alpha,min} \psi,d^1_{\alpha,min} \varphi)$. Then, $F\in \Dom (T_{\alpha,min}).$

\emph{Second inclusion.}
Let $F=(f,\varphi,\psi)\in \Dom (T_{\alpha,min})$. It means that
there exists a sequence
$((f_n,\varphi_n,\psi_n))_n $ in $\Cc^c(\Vc)\times\Cc^c(\Ec)\times\Cc^c(\Fc)$
with
\begin{equation*}
  \lim_{n\to\infty}\|f_n-f\|_{l^2(\Vc)}+\|\varphi_n-\varphi\|_{l^2(\Ec)}+\|\psi_n-\psi\|_{l^2(\Fc)}=0
\end{equation*}
and
\begin{equation*}
  \lim_{n\to\infty}\|\delta^0_\alpha(\varphi_n-\varphi)\|_{l^2(\Vc)}=0,\;\lim_{n\to\infty}\|d^1_\alpha(\varphi_n-\varphi)\|_{l^2(\Fc)}=0,\;
  \lim_{n\to\infty}\|d^0_\alpha(f_n-f)+\delta^1_\alpha(\psi_n-\psi)\|_{l^2(\Ec)}=0.
\end{equation*}
We obtain directly that
$\varphi\in \Dom (\delta^{0}_{\alpha,min})\cap  \Dom (d^{1}_{\alpha,min})$.
  \\
  We have to show now that $d^0_\alpha f_n$ and $\delta^1_\alpha\psi_n$ are Cauchy
  sequences.
\begin{Lemm}In the situation of the theorem, there exists a positive constant $C'$
  such that for any $g\in\Cc^c(\Vc)$ and $\eta\in\Cc^c(\Fc)$
  \begin{equation}\label{lem}
   |\seq{d^0_\alpha g,\delta^1_\alpha
     \eta}_{l^2(\Ec)}|\leq C'\|g\|_{l^2(\Vc)}\,\|\eta\|_{l^2(\Fc)}.
 \end{equation}
\end{Lemm}
\begin{Demo}. Let $g\in\Cc^c(\Vc)$ and $\eta\in\Cc^c(\Fc)$, then
  \begin{multline*}
    \seq{d^0_\alpha g,\delta^1_\alpha \eta}_{l^2(\Ec)}=\seq{d^1_\alpha\circ d^0_\alpha g,\eta}_{l^2(\Fc)}\\
    =-\sum_{(x,y,z)\in\Fc}s(x,y,z)\bar\eta(x,y,z)\sin\left(\frac{\widehat{\alpha}(x,y,z)}{6}\right)
    \hspace{2.4cm}\\
     \left(\e^{\frac{i}{3}(\alpha(x,z)+\alpha(x,y))}g(x)
     +\e^{\frac{i}{3}(\alpha(y,z)+\alpha(y,x))}g(y)+
     \e^{\frac{i}{3}(\alpha(z,x)+\alpha(z,y))}g(z)\right)
\end{multline*}
by the calculus of section \ref{magfield}. Now using the hypothesis of bounded
curvature (\ref{eq:bc}) one has
\begin{equation*}
  \begin{split}
\sum_{(x,y,z)\in\Fc}s(x,y,z)&\sin^2\left(\frac{\widehat{\alpha}(x,y,z)}{6}\right)
  \left(|g(x)|^2+|g(y)|^2+|g(z)|^2\right)\\
  &=3\sum_{x\in\Vc}|g(x)|^2\sum_{e\in\Fc_{x}}s(x,e)\sin^2\left(\frac{\widehat{\alpha}(x,y,z)}{6}\right)\\
    &\leq 3C\|g\|^2_{l^2(\Vc)}.
   \end{split}
 \end{equation*}
 The Cauchy-Schwarz inequality gives finally the result with $C'=3\sqrt{C}$.
\end{Demo}

 Now we return to the proof of the theorem: let $m,n\in \N$
 \begin{align*}
   \|d^0_\alpha(f_n-f_m)\|^2_{l^2(\Ec)}+&\|\delta^1_\alpha(\psi_n-\psi_m)\|^2_{l^2(\Ec)}\\
   =\|d^0_\alpha(f_n&-f_m)+\delta^1_\alpha(\psi_n-\psi_m)\|^2_{l^2(\Ec)}-2\seq{d^0_\alpha(f_n-f_m),\delta^1_\alpha(\psi_n-\psi_m)}_{l^2(\Ec)}\\
   \leq &\|d^0_\alpha(f_n-f_m)+\delta^1_\alpha(\psi_n-\psi_m)\|^2_{l^2(\Ec)}+2C'\|f_n-f_m\|_{l^2(\Vc)}\,\|\psi_n-\psi_m\|_{l^2(\Fc)}.
 \end{align*}

 We conclude that both $(d^0_\alpha f_n)_n$ and $(\delta_1\psi_n)_n$ converge so
 \begin{equation*}
   f\in \Dom (d^{0}_{\alpha,min})\quad \hbox{and}\quad\psi\in  \Dom (\delta^{1}_{\alpha,min}).
 \end{equation*}

\end{Demo}

\section{Application to magnetic triangulations}
The aim of this section is to give a concrete way to prove that the notion of $\chi-$completeness
covers the $\chi_{\alpha}$-completeness studied in \cite{ABDE}. Let $\alpha$ be an
arbitrary magnetic potential on a given triangulation $\Tc_{\alpha}$ a and fix an origin vertex
$\Oc\in\Vc$.
The combinatorial distance $d_{\rm comb}$ between two vertices $x,~y$ is given by

\[d_{\rm comb}(x,y):=\displaystyle\min_{\gamma\in\Gamma_{xy}}\mathcal{L}(\gamma)
\]
where $\Gamma_{xy}$ is the set of all paths from $x$ to $y$ and $\mathcal{L}(\gamma)$ denotes
the length (the number of edges) of a path $\gamma.$ To simplify the notation, we write
$d_{\rm comb}(\Oc,x)$ by $|x|.$
We denote by $\Bc_{n}$ the ball with center the origin vertex $\Oc$ and radius $n$, i.e
$\Bc_{n}:=\left\{x\in\Vc:~|x|\leq n\right\}.$
\begin{theo}\label{thh}
  Let $\Tc_{\alpha}$ be a weighted magnetic triangulation endowed with an origin. Assume that
  $$\displaystyle\sup_{x\in\Bc_{n}}deg_{\Vc}(x)=O(n^{2})\hbox{ and }
\displaystyle\sup_{(x,y)\in\Bc_{n}\times\Bc_{n\pm1}}deg_{\Ec}(x,y)=O(n^{2})
  \mbox{ when } n\longrightarrow\infty.
  $$
  Then $\Tc_{\alpha}$ is $\chi-$complete.
\end{theo}
\begin{Demo}

Let $n\in\N,$ we define the cut-off function $\chi_{n}:\Vc\longrightarrow\R$ by
\[\chi_{n}(x)
:=\left((2-\dfrac{|x|}{n+1})\vee0\right)\wedge1.
\]

If $x\in\Bc_{n+1},$ we have that
$\chi_{n}(x)=1$ and if $x\in\Bc^{c}_{2(n+1)},$ we have that $\chi_{n}(x)=0.$ Then, $\chi_{n}$ has a finite support.

We remark that if $y\sim x$ then, by the triangle inequality, $||x|-|y||\leq 1$ so
\begin {equation}\label{eq:dchi}
  y\sim x\Rightarrow |d^{0}\chi_{n}(x,y)|=|\chi_{n}(x)-\chi_{n}(y)|\leq\frac{1}{n+1}.
\end {equation}

We conclude that, thanks to the hypothesis on the degrees, there exists a constant $C$
such that for any $x\in\Vc$
\begin{equation*}
  \dfrac{1}{c(x)}\displaystyle\sum_{y\sim x}r(x,y)\left|d^{0}\chi_{n}(x,y)\right|^{2}
  \leq\dfrac{deg_{\Vc}(x)}{(n+1)^{2}}\leq C.
\end{equation*}

Now, let $(x,y)\in\Ec.$ Using again \eref{eq:dchi},
we have for some $M>0$ independent from $(x,y)$ and $n$

\begin{equation*}
  \begin{split}
\dfrac{1}{r(x,y)}\displaystyle\sum_{t\in\Fc_{xy}}s(x,y,t)
|d^{0}\chi_{n}(t,x)-&d^{0}\chi_{n}(t,y)|^{2}\\
& \leq \dfrac{2}{r(x,y)}\displaystyle\sum_{t\in\Fc_{xy}}s(x,y,t)\left(\left|d^{0}\chi_{n}(t,x)\right|^{2}
+\left|d^{0}\chi_{n}(t,x)\right|^{2}\right)\\
& \quad= \dfrac{4}{r(x,y)}\displaystyle\sum_{t\in\Fc_{xy}}s(x,y,t)\left|\chi_{n}(x)-\chi_{n}(t)\right|^{2}\\
& \leq\dfrac{4\, deg_{\Ec}(x,y)}{(n+1)^{2}}\\
& \leq M.
  \end{split}
\end{equation*}

The definition of $\chi-$completeness is satisfied.
\end{Demo}

\subsection{Book-like triangulations}
We recall the definition of 1-dimensional decomposition given in \cite{BG} for graphs.
This notion generalizes the decomposition in \emph{spheres} on a graph with origin, where the spheres
are defined by
$$\Sc_{n}:=\left\{x\in\Vc:~|x|=n\right\}.
$$

\begin{defi}
A 1-dimensional decomposition of a graph $\left(\Vc,\Ec\right)$ is a family of finite sets $\left(\Sc_{n}\right)_{n\in\N}$ which forms a partition of $\Vc$ and such that
for all $(x,y)\in\Sc_{n}\times\Sc_{m},$
$$(x,y)\in\Ec\Longrightarrow|n-m|\leq1.
$$
\end{defi}

The following definition is introduced in \cite{Che}.
\begin{defi}
  We say that a triangulation is a book-like triangulation endowed with an origin $\Oc$ (Fig.2),
  if there exists
  a 1-dimensional decomposition $\left(\Sc_{n}\right)_{n\in\N}$ of its graph such that for all
  $n\in\N,$
\begin{itemize}
  \item[\emph{i})]$\Sc_{0}=\left\{\Oc\right\},$~~ $\left|\Sc_{2n+1}\right|=2$
  and if $x,y\in\Sc_{2n+1}$ and $x\neq y$ then $(x,y)\in\Ec$,
  \item[\emph{ii})]$\forall x,y\in\Sc_{2n+2},\;(x,y)\notin\Ec$,
  \item[\emph{iii})]$\forall x,y\in\Sc_{2n+1},\,x\neq y$ and $\forall z\in\Sc_{2n}\cup\Sc_{2n+2}$
then $(x,y,z)\in\Fc$.
\end{itemize}
\end{defi}

\begin{figure}[!t]
\centering
\begin{minipage}[t]{10cm}
\includegraphics*[width=9cm,height=4.5cm]{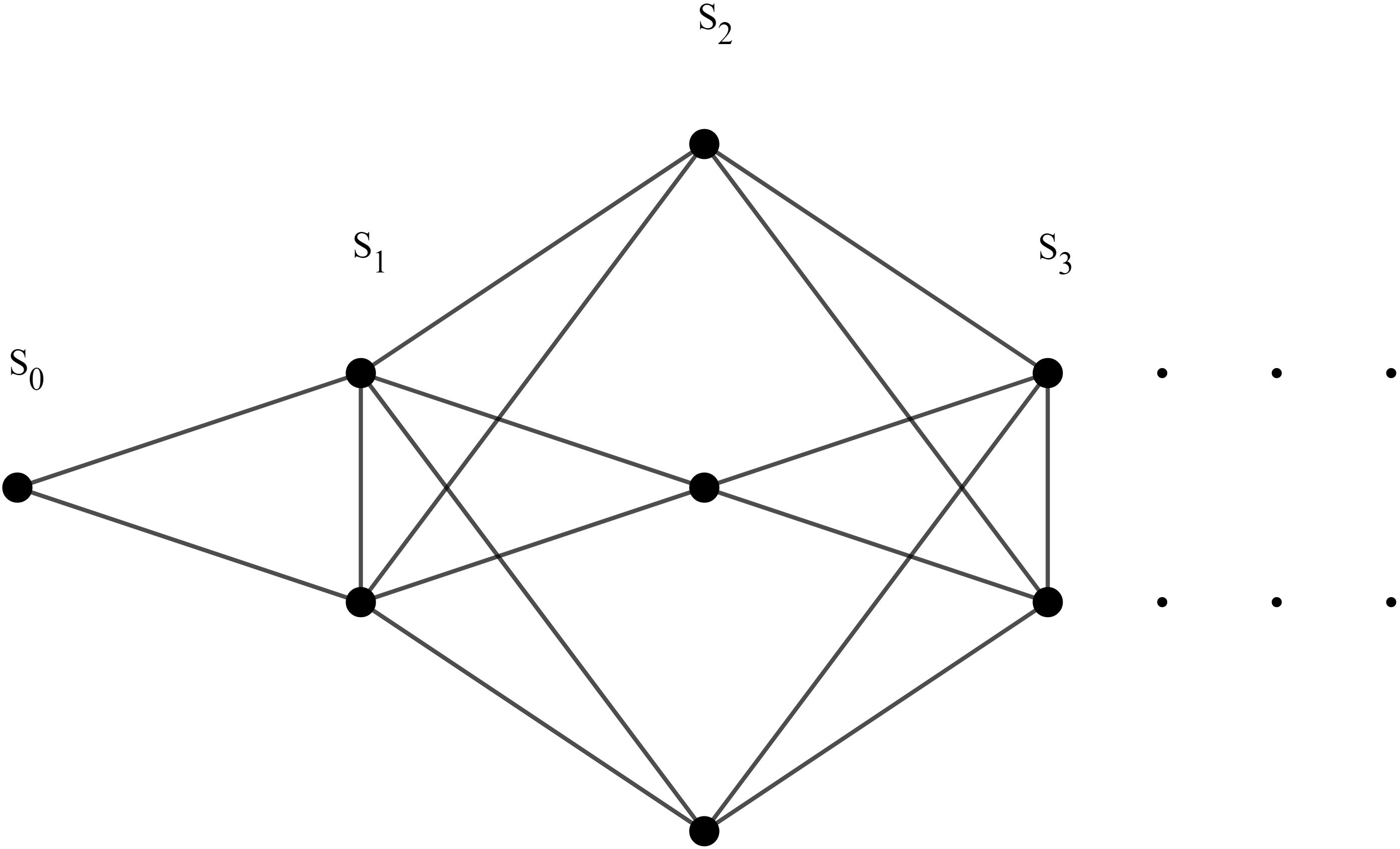}
\caption{A book-like triangulation}
\end{minipage}
\end{figure}

\begin{Exp}
Let $\Tc_{\alpha}$ be a magnetic book-like triangulation endowed with an origin $\Oc$
and an arbitrary magnetic potential $\alpha.$ Take $0<\beta\leq2.$ We set $c(x)=1$ and
$val(x):=\sharp\left\{y\in\Vc,~y\sim x\right\},$ for all $x\in\Vc.$ If $\Tc_{\alpha}$
satisfies
\begin{align*}val(x)&=\lfloor (2n+1)^{\beta}\rfloor+4,
  \mbox{ for all }n\in\N\mbox{ and }x\in\Sc_{2n+1}\\
r(x,y)&=\dfrac{val(x)val(y)}{|x||y|},
~\mbox{ for all }(x,y)\in\Ec\\
\hbox{and }\hspace{2cm}
s(x,y,z)&=r(x,y)r(y,z)r(z,x),\mbox{ for all }(x,y,z)\in\Fc.
\end{align*}

Then an application of Theorem \ref{thh} gives that $\Tc_{\alpha}$ is $\chi$-complete.\\

Indeed,  we remark first that $val(x)=4$ for $x\in\Sc_{2n}, n\geq 1)$ and 
\begin{equation}\label{exe}
val(x)=|\Sc_{2n}|+|\Sc_{2n+2}|+1,\mbox{ for all }n\in\N
\mbox{ and for }x\in\Sc_{2n+1}.
\end{equation}

Let $n\in\N^{*}$, the equation (\ref{exe}) gives:

\emph{First case.} If $x\in\Sc_{2n+1},$ let $\bar x$ be the other vertex in $\Sc_{2n+1}$.
There is some $C>0$ independent from $x$ and $n$ such that
\begin{equation*}
  \begin{split}
    \displaystyle\sum_{y\sim x}\dfrac{r(x,y)}{c(x)}
    & = \left(\dfrac{val(x)}{|x|}\right)\left(\dfrac{val(\bar x)}{|\bar x|}\right)+
    \dfrac{val(x)}{|x|}\displaystyle\sum_{y\in\Sc_{2n}\cup\Sc_{2n+2}}\dfrac{val(y)}{|y|}\\
    & = \left(\dfrac{\lfloor (2n+1)^{\beta}\rfloor+4}{2n+1}\right)^{2}
    +\dfrac{\lfloor (2n+1)^{\beta}\rfloor+4}{2n+1}\left(\dfrac{4|\Sc_{2n}|}{2n}+
    \dfrac{4|\Sc_{2n+2}|}{2n+2}\right)\\
    & \leq \left(\dfrac{\lfloor (2n+1)^{\beta}\rfloor+4}{2n+1}\right)^{2}
    +2\dfrac{\lfloor (2n+1)^{\beta}\rfloor+4}{2n+1}\left(\dfrac{|\Sc_{2n}|+|\Sc_{2n+2}|}{n}\right)\\
    & \leq \left(\dfrac{\lfloor (2n+1)^{\beta}\rfloor+4}{2n+1}\right)^{2}
    +2\dfrac{\left(\lfloor (2n+1)^{\beta}\rfloor+4\right)^{2}}{n(2n+1)}\\
    & \leq Cn^{2\beta- 2}  \hspace{1cm}   (  0<\beta\leq2  )  \\
    &\leq C n^{2}.
  \end{split}
\end{equation*}

\emph{Second case.} If $x\in\Sc_{2n},$ we have for some $C'>0$ independent from $x$ and $n$
\begin{equation*}
  \begin{split}
    \displaystyle\sum_{y\sim x}\dfrac{r(x,y)}{c(x)}
    & =\dfrac{val(x)}{|x|}\displaystyle\sum_{y\sim x,~y\in\Sc_{2n\pm1}}\dfrac{val(y)}{|y|}\\
    & =\dfrac{4}{2n}\left(\dfrac{|\Sc_{2n-1}|\left(\lfloor (2n-1)^{\beta}\rfloor+4\right)}{2n-1}
    +\dfrac{|\Sc_{2n+1}|\left(\lfloor (2n+1)^{\beta}\rfloor+4\right)}{2n+1}\right)\\
    & =\dfrac{4}{n}\left(\dfrac{\lfloor (2n-1)^{\beta}\rfloor+4}{2n-1}
    +\dfrac{\lfloor (2n+1)^{\beta}\rfloor+4}{2n+1}\right)\\
    & \leq C'n^{\beta- 2}\leq C',   \; \;\;  \; \;\;   (  0<\beta\leq2  )  .
  \end{split}
\end{equation*}

On the other hand, let $(x,y)\in\Sc^{2}_{2n+1},\,x\neq y$. We have for some $M>0$
independent from $(x,y)$ and $n$
\begin{equation*}
  \begin{split}
    \displaystyle\sum_{t\in\Fc_{xy}}\dfrac{s(x,y,t)}{r(x,y)}
    & = \displaystyle\sum_{t\in\Fc_{xy}}r(y,t)r(t,x)\\
    & = \dfrac{val(x)val(y)}{|x||y|}\displaystyle\sum_{t\in\Sc_{2n}\cup\Sc_{2n+2}}
    \left(\dfrac{val(t)}{|t|}\right)^{2}\\
    & = \left(\dfrac{\lfloor (2n+1)^{\beta}\rfloor+4}{2n+1}\right)^{2}
    \left(\dfrac{16|\Sc_{2n}|}{(2n)^2}+\dfrac{16|\Sc_{2n+2}|}{(2n+2)^2}\right)\\
    & \leq 16\left(\dfrac{\lfloor (2n+1)^{\beta}\rfloor+4}{2n+1}\right)^{2}
    \left(\dfrac{|\Sc_{2n}|+ |\Sc_{2n+2}|}{(2n)^2}\right)\\
    & \leq 4\left(\dfrac{\lfloor (2n+1)^{\beta}\rfloor+4}{2n+1}\right)^{2}
    \left(\dfrac{\lfloor (2n+1)^{\beta}\rfloor+4}{n^2}\right)\\
    & \leq M n^{3\beta- 4},   \; \;\;  \; \;\;   (  0<\beta\leq2  )  \\
    &\leq M n^{2}.
  \end{split}
\end{equation*}

If $(x,y)\in\Sc_{2n}\times\Sc_{2n+1},$ we have for some $M'>0$ independent from $(x,y)$ and $n$
\begin{equation*}
  \begin{split}
    \displaystyle\sum_{t\in\Fc_{xy}}\dfrac{s(x,y,t)}{r(x,y)}
    & = \displaystyle\sum_{t\in\Fc_{xy}}r(y,t)r(t,x)\\
    & = \dfrac{val(x)val(y)}{|x||y|}\displaystyle\sum_{t\in\Sc_{2n+1}}
    \left(\dfrac{val(t)}{|t|}\right)^{2}\\
    & = \dfrac{4
    \left(\lfloor (2n+1)^{\beta}\rfloor+4\right)}{2n(2n+1)}
   2 \left(\dfrac{\lfloor (2n+1)^{\beta}\rfloor+4}{2n+1}\right)^{2}\\
    & \leq M' n^{3\beta- 4}\hspace{1cm}   (  0<\beta\leq2  )             \\
    &\leq M' n^{2}.
  \end{split}
\end{equation*}
\end{Exp}

\subsection{Case of a 1-dimensional decomposition triangulation}
In this section under a specific choice of magnetic potential $\alpha,$ we construct
an example of a triangulation that is $\chi$-complete and not $\chi_{\alpha}$-complete.
The following definition is introduced in \cite{Che}.
\begin{defi}
We say that a triangulation is a \emph{1-dimensional decomposition triangulation} if there exists a
1-dimensional decomposition $\left(\Sc_{n}\right)_{n\in\N}$ of its graph (Fig.3).
\end{defi}

We now divide the degree with respect to the simple 1-dimensional decomposition triangulation
$$\eta_n^{\pm}:=\displaystyle\sup_{x\in\Sc_n}\left|\Vc(x)\cap\Sc_{n\pm1}\right|,
~\beta_n^{\pm}:=\sup_{e\in\Sc_n\times\Sc_{n^{\pm}1}}\left|\Fc_e\right|,
~\gamma_n^{\pm}:=\sup_{e\in\Sc_n^2}\left|\left(\left(\Fc_e\cap{\Sc_{n\pm1}}\right)\times\Sc_{n}^2\right)\cap\Fc\right|.
$$

\begin{figure}[!t]
\centering
\begin{minipage}[t]{10cm}
\includegraphics*[width=11cm,height=7cm]{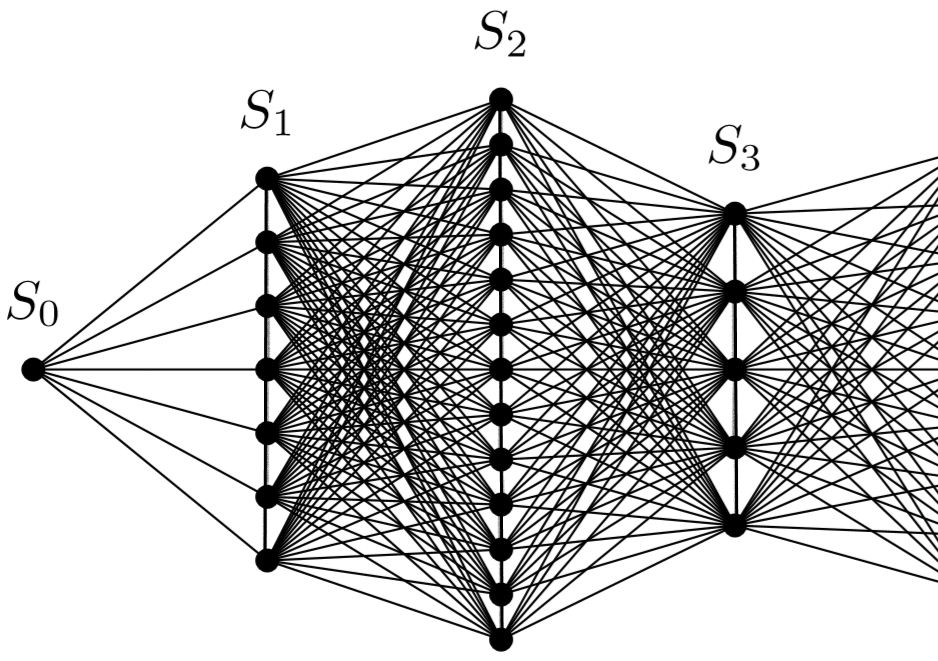}
\caption{A 1-dimensional decomposition triangulation}
\end{minipage}
\end{figure}

\begin{Exp}
Let $\Tc_{\alpha}$ be a simple 1-dimensional decomposition magnetic triangulation and
$$\alpha_{xy}:=\left(|x|-|y|\right)\pi.
$$
\end{Exp}
Thus,
$$\displaystyle\sum_{y\sim x}\sin^{2}\left(\frac{\alpha_{xy}}{2}\right)=val(x)
$$

for all $x\in\Vc.$ Using (\ref{eq:1}), if $val(.)$ is unbounded,
then $\Tc_{\alpha}$ is not $\chi_{\alpha}$-complete.

In contrast, in  \cite[Thm 6.2]{Che}, the author proves that if
$$\displaystyle\sum_{n\in\N}\dfrac{1}{\sqrt{\xi(n,n+1)}}=\infty
$$

where $\xi(n,n+1)=\eta_n^{+}+\eta_{n+1}^{-}+\beta_n^{+}+\gamma_n^{+}+\gamma_{n+1}^{-}$,
then $\Tc_{\alpha}$ is $\chi-$complete.


As a consequence, the operators $T_\alpha$ and $\Delta_\alpha$ are essentially self-adjoint on
$\Cc^c(\Vc)\oplus\Cc^c(\Ec)\oplus\Cc^c(\Fc)$.
But we remark also that the magnetic field is trivial, it is
the differential of a function, then its holonomy is null and the curvature also. We
can then apply Theorem \ref{min-bc} to $T_\alpha$.

\medskip
\textbf{\textit{\emph{Acknowledgments}}:}
The authors would like to thank \textit{Luc Hillairet} for fruitful discussions and helpful comments.
And we sincerely thank the \textit{anonymous referee} for both the careful reading and the useful guidance as well as the words of encouragement on our paper.

\end{document}